\documentclass[11pt]{article}
\usepackage{amssymb}
\usepackage{amsmath}
\usepackage{amsthm}
\newcommand{\cA}{{\cal A}}

\newcommand{\cC}{{\cal C}}
\newcommand{\cD}{{\cal D}}

\newcommand{\cO}{{\cal O}}
\newcommand{\cL}{{\cal L}}
\newcommand{\cM}{{\cal M}}
\newcommand{\cN}{{\cal N}}
\newcommand{\cF}{{\cal F}}
\newcommand{\cK}{{\cal K}}
\newcommand{\cP}{{\cal P}}
\newcommand{\cQ}{{\cal Q}}
\newcommand{\cR}{{\cal R}}
\newcommand{\cS}{{\cal S}}
\newcommand{\cT}{{\cal T}}

\newcommand{\cX}{{\cal X}}
\newcommand{\cY}{{\cal Y}}

\newcommand{\EE}{{\mathbb E}}
\newcommand{\NN}{{\mathbb N}}
\newcommand{\ZZ}{{\mathbb Z}}

\newcommand{\CC}{{\mathbb C}}

\newcommand{\DD}{{\mathbb D}} 

\newcommand{\gm}{\mathfrak{m}}    
\newcommand{\gn}{\mathfrak{n}}
\newcommand{\gp}{\mathfrak{p}}
\newcommand{\gq}{\mathfrak{q}}

\newcommand{\gr}{\mathfrak{r}}

\newcommand{\on}{\operatorname}
\newcommand{\Sh}{\on{Sh}}
\newcommand{\Fr}{\on{Fr}}

\newcommand{\Fl}{{\cal F}l}

\newcommand{\Qlb}{\mathbb{\bar Q}_\ell}
\newcommand{\Gm}{\mathbb{G}_m}
\newcommand{\A}{\mathbb{A}}
\newcommand{\Ql}{\mathbb{Q}_\ell}
\newcommand{\toup}[1]{\stackrel{#1}{\to}}
\newcommand{\hook}[1]{\stackrel{#1}{\hookrightarrow}}

\newcommand{\getsup}[1]{\stackrel{#1}{\gets}}

\newcommand{\Hom}{\on{Hom}}

\newcommand{\Mod}{\on{Mod}}
\renewcommand{\mod}{\on{mod}}
\newcommand{\Ext}{\on{Ext}}
\newcommand{\Flag}{\on{Flag}}
\newcommand{\End}{\on{End}}
\newcommand{\Ann}{\on{Ann}}
\newcommand{\Sym}{\on{Sym}}

\newcommand{\mult}{\on{mult}}
\newcommand{\Aut}{\on{Aut}}
\newcommand{\Autb}{\on{\underline{Aut}} }
\newcommand{\RG}{\on{R\Gamma}}

\newcommand{\Perv}{\on{Perv}}

\newcommand{\parf}{\on{parf}}

\newcommand{\sym}{\on{sym}}

\newcommand{\qcoh}{\on{qcoh}}
\renewcommand{\part}{\on{part}}
\newcommand{\Spr}{{{\cal S}pr}}

\newcommand{\Pic}{\on{Pic}}
\newcommand{\uPic}{\on{\underline{Pic}}}
\newcommand{\Bun}{\on{Bun}}
\newcommand{\Bunb}{{\on{\underline{Bun}} }}
\newcommand{\multb}{{\on{\underline{mult}} }}
\newcommand{\rk}{\on{rk}}
\newcommand{\Spec}{\on{Spec}}
\newcommand{\Specf}{\on{Spf}}

\newcommand{\supp}{\on{supp}}

\newcommand{\HOM}{{{\cal H}om}}
\newcommand{\END}{{{\cal E}nd}}

\newcommand{\Ob}{\on{Ob}}
\newcommand{\Sets}{\on{Sets}}
\newcommand{\GL}{\on{GL}}
\newcommand{\PGL}{\on{PGL}}

\newcommand{\pr}{\on{pr}}
\newcommand{\id}{\on{id}}

\newcommand{\tr}{\on{tr}}

\newcommand{\QED}{$\square$} 
\newcommand{\Fq}{\mathbb{F}_q}  
\newcommand{\Fp}{\mathbb{F}_p}  
\newcommand{\iso}{{\widetilde\to}}

\newcommand{\comp}{\circ}

\renewcommand{\H}{\on{H}}   
\newcommand{\R}{\on{R}\!}   
\newcommand{\Lotimes}{\stackrel{L}{\otimes}}
\newcommand{\Cotimes}{\mathop{{\hat\otimes}}}
\newcommand{\D}{\on{D}}       
\newcommand{\wt}{\widetilde}
\newcommand{\select}[1]{{\it{#1}}}
\newcommand{\p}{\prime}
\renewcommand{\u}[1]{\underline{#1}}
\newcommand{\otimesl}{\mathop{\otimes}\limits}

\newcommand{\Isom}{\on{\underline{Isom}}}
\newcommand{\V}{\mathbb{V}}

\renewcommand{\div}{\on{div}}
\renewcommand{\P}{\mathbb{P}}

\newtheorem{Lm}{Lemma}

\newtheorem{MGTh}{Main Global Theorem{\hskip -3pt}}

\newtheorem{Pp}{Proposition}
\newtheorem{Cor}{Corolary}
\newtheorem{Con}{Conjecture}
\newtheorem{Slm}{Sublemma}

\newtheorem{lm}{Lemma}

\newtheorem{Lm_section}{Lemma}[section]

\theoremstyle{remark}
\newtheorem{Rem}{Remark}
\newtheorem{Rems}{Remarks}

\theoremstyle{definition}
\newtheorem{Def}{Definition}

\newenvironment{Prf}{\par\noindent {\it Proof }}{\QED}
\newcommand{\Step}[1]{\par\noindent{\bf Step {#1}}.}

\begin{document}
\author{Sergey Lysenko}
\title{Global geometrised Rankin-Selberg method for $\GL(n)$}
\date{}
\maketitle
\begin{abstract}
\noindent {\scshape Abstract} \hskip 0.8 em 
We propose a geometric interpretation of the classical Rankin-Selberg 
method for $\GL(n)$ in the framework of the geometric Langlands
program. We show that the geometric Langlands conjecture for
an irreducible unramified local system $E$ of rank $n$ on a curve 
implies the existence of automorphic sheaves corresponding to the 
universal deformation of $E$. Then we calculate the `scalar product' 
of two automorphic sheaves attached to this universal deformation.
\end{abstract}

\section*{Introduction}
0.1. In this paper, which is a continuation of \cite{Ly2}, 
we propose a geometric interpretation of the global Rankin-Selberg method
for $\GL(n)$ in the framework of the geometric Langlands program. 
The corresponding classical result computes the scalar product of two
cuspidal nonramified automorphic forms on $\GL(n)$ over a function field. 

 Our first motivation is the following result of G. Laumon (\cite{L4}) and
M. Rothstein (\cite{R}) in the case of $\GL(1)$. Let $X$ be a smooth, projective,
connected curve of genus $g\ge 1$ over $\CC$. Let $M'$ be the Picard scheme 
classifying invertible $\cO_X$-modules $L$ of degree zero. Denote by $M$ the
coarse moduli space of invertible $\cO_X$-modules $L$ with connection
$\nabla:L\to L\otimes_{\cO_X}\Omega_X$. This is an abelian group scheme over
$\CC$ (for the tensor product), which has a natural structure of a
$\H^0(X,\Omega_X)$-torsor over $M'$. 
In \cite{L4,R} a certain invertible $\cO_{M\times
M'}$-module $\Aut$ with connection (relative to $M$) is considered as a
kernel of two Fourier transforms $\cF: \D^b_{\qcoh}(\cD_{M'})\to \D^b_{\qcoh}(\cO_M)$
and $\cF': \D^b_{\qcoh}(\cO_M)\to \D^b_{\qcoh}(\cD_{M'})$. 
The theorem of G. Laumon and M. Rothstein claims that \select{these 
functors are quasi-inverse to each other}. 

This result can be obtained as a formal consequence of two
orthogonality relations. One of them states that 
\select{the complex
$$
\R(\pr_{12})_*(\pr^*_{13}\Aut\otimes\pr^*_{23}\Aut)
$$
is canonically isomorphic to $\vartriangle_*\!\cO_M$ in $\D^b_{\qcoh}(\cO_{M\times
M})$ (up to a shift and a sign)}, where 
$\pr_{13},\pr_{23}:M\times M\times M'\to M\times M'$ and $\pr_{12}: M\times
M\times M'\to M\times M$ are the projections, $\vartriangle:M\to M\times M$ is
the diagonal, and the functor $\R(\pr_{12})_*$ is understood in the $\cD$-modules
sense.  

\medskip
\noindent
0.2. An $\ell$-adic analogue of this orthogonality relation is the following. Let $X$
be a smooth, projective, connected curve of genus $g\ge 1$ over an algebraically closed
field $k$. Fix a prime $\ell$ invertible in $k$ and 
an algebraic closure $\Qlb$ of $\Ql$. 
Let $E_0$ be a smooth $\Qlb$-sheaf on $X$ of rank 1. 

 The moduli space of $\ell$-adic local systems on $X$ is not known to
exist. However, we can consider the deformations of $E_0$ over $\Qlb$ (cf. Sect.
3.1). The local system $E_0$ admits a unversal deformation $E$ (cf.
Proposition~\ref{pro-repr}). Let $\Specf(R)$ be the base of this deformation. In fact,
$R$ is isomorphic to the ring of formal power series over $\Qlb$ of dimension $2g$.  

 For a positive integer $d$ we have the $R$-sheaf $E^{(d)}$ on 
the symmetric power $X^{(d)}$ of $X$ (cf. Sect. 1.5). Let $\uPic^d X$
denote the Picard scheme of $X$ parametrizing the isomorphism classes of invertible
$\cO_X$-modules of degree $d$. According to the geometric abelian class field theory,
$E^{(d)}$ descends to a smooth $R$-sheaf $E^d$ of rank 1 on $\uPic^d X$. Denote by
$E^d_1, E^d_2$ the two liftings of $E^d$ to $\Specf(R\Cotimes_{\Qlb} R)$. Essentially,
we show that \select{there is a canonical isomorphism
of $R\Cotimes R$-modules
$$
\H^{2g}(\uPic^d X, \HOM(E^d_1,E^d_2))\iso R(-g),
$$
where the $R\Cotimes R$-module structure on $R$ is given by the diagonal map $R\Cotimes
R\to R$. Besides, $\H^i(\uPic^d X, \HOM(E^d_1,E^d_2))=0$ for $i\ne 2g$.} Actually, a
bit different statement is proved (cf. the discussion at the end of Sect.~0.4).

 This is the particular case $n=1$ of our Main Global Theorem (cf.
Sect. 4.1), which is an analogue of this orthogonality relation for
$\GL(n)$.  

 Already in the case of $\GL(1)$ we observe the important role of deformations. Applying
the base change theorem for the above result, we get the scalar square of
$E^d_0$, namely
$$
\RG(\uPic^d X, \END(E_0^d,E_0^d))\iso R\Lotimes_{R\Cotimes R}\Qlb (-g)[-2g]
$$ 
As is easy to see, this complex has cohomology groups in all degrees $0,1,\ldots, 2g$.  
For $\GL(n)$ the answer
will also be simplier when considered as an object on the moduli of parameters. 

\medskip\noindent
0.3. Let $X$ be a smooth, projective, geometrically connected
curve over $\Fq$. Let $\ell$ be a prime invertible in $\Fq$. According to the
Langlands correspondence for $\GL(n)$ over function fields (proved in full generality by
L. Lafforgue), to any smooth geometrically irreducible $\Qlb$-sheaf $E$ 
on $X$ is associated a (unique up to a multiple)
cuspidal automorphic form $\varphi_E: \Bun_n(\Fq)\to \Qlb$, which
is a Hecke eigenvector with respect to $E$. The function $\varphi_E$ is defined on the
set $\Bun_n(\Fq)$ of isomorphism classes of rank $n$ vector bundles on $X$. 

 The global Rankin-Selberg for $\GL(n)$ allows to
calculate for any integer $d$ the scalar product of two (appropriately normalized)
automorphic forms
\begin{equation}
\label{scalar_square}
\sum_{L\in\Bun^d_n(\Fq)} \frac{1}{\#\Aut L}\; \varphi_{E_1^*}(L)\varphi_{E_2}(L),
\end{equation}
where $\Bun^d_n(\Fq)$ is the set of isomorphism classes of vector bundles $L$ on $X$
of rank $n$ and degree $d+n(n-1)(g-1)$, and $\#\Aut L$ stands for the number of 
elements in $\Aut L$. More precisely, this scalar product vanishes if and only if 
$E_1$ and $E_2$ are non isomorphic. In the case $E_1\iso E_2\iso E$ the answer is 
expressed in terms of the action of the geometric Frobenius endomorphism on 
$\H^1(X\otimes\bar\Fq,\,\END E)$.

Assume that $\varphi_E$ is canonically normalized (cf. footnote 1).
The local part of the classical Rankin-Selberg method for $\GL(n)$ may be stated as
the equality of formal series (cf. \cite{Ly2})
\begin{equation}
\label{equality_formal_series}
\sum_{d\ge 0}\;\; \sum_{(\Omega^{n-1}\hook{}L)\in
\,{_n\cM_d}(\Fq)}\;\;\frac{1}{\#\Aut(\Omega^{n-1}\hook{}L)}
\varphi_{E_1^*}(L)\varphi_{E_2}(L)t^d=q^{n^2(1-g)}\,
L(E_1^*\otimes E_2,t)
\end{equation}
Here $_n\cM_d(\Fq)$ is the set of isomorphism classes of pairs
$(\Omega^{n-1}\hook{}L)$, where $L$ is a vector bundle on $X$ of rank $n$ and
degree $d+n(n-1)(g-1)$, and $\Omega$ is the
canonical invertible sheaf on $X$ ($\Omega^{n-1}$ is embedded in $L$ as a subsheaf, 
i.e., the quotient is allowed to have torsion). We have denoted by
$$
L(E_1^*\otimes E_2,t)=\sum_{d\ge 0}\sum_{D\in X^{(d)}(\Fq)} 
\tr(\Fr, (E_1^*\otimes E_2)^{(d)}_D)t^d
$$
the L-function attached to the local system $E_1^*\otimes E_2$ (here $\Fr$ is the
geometric Frobenius endomorphism).  

 The calculation of (\ref{scalar_square}) is an
asymptotic argument. First, rewrite (\ref{equality_formal_series}) as
\begin{multline*}
\sum_{d\ge 0}\; \sum_{L\in \Bun_n^d(\Fq)}\;
\frac{1}{\#\Aut L}(q^{\dim\Hom(\Omega^{n-1},L)}-1)
\varphi_{E_1^*}(L)\varphi_{E_2}(L)t^d=\\ 
q^{n^2(1-g)}\, L(E_1^*\otimes E_2,t) 
\end{multline*}
The cuspidality of $\varphi_E$ implies that if $\varphi_E(L)\ne 0$ and $\deg L$ is
large enough then $\Ext^1(\Omega^{n-1},L)=0$, and 
$\dim\Hom(\Omega^{n-1},L)=d-n^2(g-1)$. To conclude, it remains to study the 
asymptotic behaviour of the above series when $t$ goes to $q^{-1}$, using the
cohomological interpretation 
$$
L(E_1^*\otimes E_2,t)=\prod_{r=0}^2\det(1-\Fr\, t,\, \H^r(X\otimes\bar \Fq, 
E_1^*\otimes E_2))^{(-1)^{r+1}}
$$
of the $L$-function.  

\medskip\noindent
0.4. Let $X$ be as in 0.2. Let $E_0$ be a smooth 
irreducible $\Qlb$-sheaf on $X$ of rank
$n$. Let $(E,R)$ be a universal deformation of $E_0$ (cf. Sect. 3.1). The base
$R$ is, in fact, isomorphic to the ring of formal power series over $\Qlb$ of
dimension $2+(2g-2)n^2$. We assume the geometric Langlands conjecture for GL(n) 
(Conjecture~\ref{Con_1}, Sect. 2.4), so that to
$E$ is associated a perverse $R$-sheaf $\Aut_E$ on the moduli stack $\Bun_n$
of rank $n$ vector bundles on $X$, which is a geometric analogue of 
the automorphic form $\varphi_E$ \footnote{We normalize $\Aut_E$ as in 
Conjecture~\ref{Con_1}, this also gives a
normalization of $\varphi_E$ as the function trace of Frobenius of $\Aut_E$.}. 
Denote by $\Aut^d_E$ the restriction of $\Aut_E$ to $\Bun^d_n$, the connected 
component of $\Bun_n$ classifying vector bundles on $X$ of rank $n$ and degree
$d+n(n-1)(g-1)$. One can calculate the cohomology
$$
\RG_c(\Bun^d_n, \,\pr_1^*\Aut^d_{E^*}\otimes_{R\Cotimes R}\pr_2^*\Aut^d_E),
$$
which is a geometric analogue of (\ref{scalar_square}). 
(Here $\pr_i: \Specf(R\Cotimes
R)\to \Specf R$ are the two projections.)
But the answer turns out to be `bad': this complex has nontrivial 
cohomology groups in infinitely many degrees
(bounded from above). To get the desired answer, we modify the problem as
follows.

 Scalar automorphisms of vector bundles provide an action of $\Gm$ on $\Bun_n$ by
2-automorphisms of stacks. We introduce the stack $\Bunb_n$ (cf.
Sect. 2.5), the quotient of $\Bun_n$ under this action. There exists
a perverse $R$-sheaf $\Autb_E$ on $\Bunb_n$ such that the inverse image of
$\Autb_E[-1](-\frac{1}{2})$ under the projection $\Bun_n\to\Bunb_n$ is 
identified with $\Aut_E$. 

Our Main Global Theorem essentially says that 
\select{if Conjecture~\ref{Con_1} is true then
for any integers $i$ and $d$ there is a canonical isomorphism 
of $R\Cotimes R$-modules
$$
\H^i_c(\Bunb^d_n, \pr_1^*\Autb^d_{E^*}\otimes_{R\Cotimes R}\pr_2^*\Autb^d_E)
\iso 
\begin{cases}
R,  & \mbox{ if } i=0\\
0,  & \mbox{ if } i\ne 0,
\end{cases}
$$
where the $R\Cotimes R$-module structure on $R$ is given by the diagonal map
$R\Cotimes R\to R$.}

 Actually a bit different statement is proved. Conceptually, we should work
in the derived categories $\D^b_c(Y,R)$ of complexes of $R$-sheaves 
in \'etale topology on $Y$, where $Y$ is a $k$-scheme of finite type
and $R$ is the ring of formal power series over $\Qlb$. (For $R=\Qlb$ this 
is the derived category of $\ell$-adic sheaves on $Y$).
However, for $\dim R>0$ the definition of this derived category is not known. 
We impose an additional condition: the reduction of $E_0$ modulo $\ell$ has no
nontrivial endomorphisms. This allows to work in another triangulated category 
denoted $\D^b_c(Y,\bar R)_{\sigma}$ (cf. Sect.~1.4), which `approximates' the 
desired one. Our Main Global Theorem is an isomorphism in such triangulated 
category.

\medskip\noindent
0.5. The paper is organized as follows. 
In Sect.~1 we partially establish some formalism of $R$-sheaves on schemes. 
We adopt the point of view that an appropriate formalism of $R$-sheaves holds 
for algebraic stacks locally of finite type over $k$. 
In Sect.~2 we formulate the geometric Langlands conjecture relative to 
$R$ (Conjecture~\ref{Con_1}), and essentially show that if this conjecture 
is true over $\Qlb$ then it is also true over the rings of formal power 
series over $\Qlb$ (actually a bit different statement is proved, 
cf. Proposition~2).  

 In Sect. 3.1-3.3 we study the deformations of local systems on $X$. 
In Sect. 3.4 we calculate the cohomology of some natural sheaves arising 
from the universal deformations. The proof of Main Global Theorem is given in 
Sect.~4.

 We refer the reader to Sect. 0.2 of \cite{Ly2} for our conventions.

\section{$R$-sheaves}

1.1 Let $\EE$ be a finite extension field of $\Ql$. Denote by
$\cC_{\,\EE}$ the category of local Artin
$\EE$-algebras with residue field $\EE$ (the morphisms are local
homomorphisms of $\EE$-algebras). 

Let $A\in \Ob(\cC_{\,\EE})$ and $Y$ be a $k$-scheme of finite type. The
category $\Sh(Y,A)$ of constructable $A$-sheaves  on $Y$ is the category of
pairs $(E,\rho)$, where $E$ is a constructable $\EE$-sheaf on $Y$, 
and $\rho: A\to \End(E)$ is an action of $A$ on $E$. 

 A constructable $A$-sheaf $(E,\rho)$ is \select{smooth of rank $m$}, 
if $E$ is a smooth $\EE$-sheaf, and the fibres of $(E,\rho)$ are free 
$A$-modules of rank $m$.   

 Let $R$ be a complete local noetherian $\EE$-algebra with
residue field $\EE$ and maximal ideal $\gm$. We let 
$\Sh(Y,R)$ be the projective 2-limit of $\Sh(Y, R/\gm^n)$. 
This is a category whose objects are projective systems 
$(F_n,\psi_n)_{n\in\NN}$, where $F_n\in \Ob(\Sh(Y,R/\gm^n))$ and 
$\psi_n: F_{n+1}\otimes_{R/\gm^{n+1}} R/\gm^n\iso F_n$ is an isomorphism.
Morphisms are defined as the morphisms of the corresponding projective
systems. 

 Conceptually, we should work in the derived category 
$\D^b_c(Y,R)$, where $R$ is the ring of formal power
series in several variables over $\Qlb$, but its definition is not known. 
We `approximate' it by another triangulated category defined below. 
 
\medskip
\noindent
1.2 Let $A$ be a commutative ring, $S\subset A$ be a multiplicative system.
We denote by $A_S$ the localization of $A$ in $S$.   Given an $A$-additive
category $\cA$, denote by
$\cA\otimes_A A_S$ the category such that $\Ob(\cA\otimes_A A_S)=\Ob(\cA)$,
and morphisms in $\cA\otimes_A A_S$ from $K$ to $K'$ are
$\Hom_{\cA}(K,K')\otimes_A A_S$. This is an $A_S$-additive category.
The following statement is left to the reader (the points i) and iii) follow
from \cite{GL}, Proposition 3.6.1). 

\begin{Lm} 
\label{Lm_for_Gabber}
1) Let $D$ be an $A$-additive triangulated category. 
Then $D\otimes_A A_S$ is triangulated (its
distinguished triangles are those isomorphic to images of distinguished
triangles in $D$). If $F:D\to D'$ is an $A$-additive triangulated functor
between $A$-additive triangulated categories then the induced functor
$D\otimes_A A_S\to D'\otimes_A A_S$ is triangulated.

\smallskip
\noindent
2) Assume, in addition, that $D$ is equiped with 
a $t$-structure $(D^{\leqslant 0}, D^{\geqslant 0})$ such that any $K\in\Ob(D)$ 
is of finite amplitude, i.e., $D=\cup_i \, D^{\leqslant i}=\cup_i 
\, D^{\geqslant i}$. Let $\cA$ denote the core of this $t$-structure. 
Let $\cS$ be the full-subcategory of $D$ consisting of $K$ with
$\H^i(K)\otimes_A A_S=0$ in $\cA\otimes_A A_S$ for all $i$. Then
\begin{itemize}
\item[i)] $\cS$ is a thick subcategory of $D$;
\item[ii)] the localization of $D$ in $\cS$ is canonically equivalent to
$D\otimes_A A_S$;
\item[iii)] $D\otimes_A A_S$ has a $t$-structure $(D^{\leqslant 0}_S, 
D^{\geqslant 0}_S)$, where $D^{\leqslant 0}_S$ (resp., $D^{\geqslant 0}_S$) 
is the essential image of those $K\in D$ that $\H^i(K)\otimes_A A_S=0$ for $i>0$ 
(resp., for $i<0$). The core of this $t$-structure is an abelian category 
equivalent to $\cA\otimes_A A_S$. \ \QED 
\end{itemize} 
\end{Lm}

\begin{Rem}
\label{Rem_strange}
We will use the following trivial observation. Let $c\!-\!A\!-\!\mod$ denote the
category of finite type $A$-modules. If $A$ is noetherian 
then the natural functor
$(c\!-\!A\!-\!\mod)\otimes_A A_S\to c\!-\!A_S\!-\mod$ is fully faithful.
\end{Rem}

\smallskip

The following statement is left to the reader (point 3) follows from
\cite{Del}, 4.3 and 4.5.1).

\begin{Lm} 
\label{Lm_perfect_derived}
Let $\Lambda$ be a local noetherian ring with
residue field $\kappa_{\Lambda}$.

\noindent
1) If $K$ is a perfect complex of $\Lambda$-modules then it can be
represented as a direct sum $K=K_1\oplus K_2$ of two perfect complexes,
where $K_1$ is acyclic, and the differential in
$K_2\otimes_{\Lambda}\kappa_{\Lambda}$ is zero.

\smallskip
\noindent
2) If $K_1$ and $K_2$ are perfect complexes of $\Lambda$-modules
such that the differential in $K_i\otimes_{\Lambda}\kappa_{\Lambda}$ is zero
($i=1,2$) and $f:K_1\to K_2$ is a homotopical equivalence then $f$ is an
isomorphism of complexes. 

\smallskip
\noindent
3) Let $\D_{parf}(\Lambda)$ be the category whose objects 
are perfect complexes of $\Lambda$-modules, and the morphisms are 
morphisms of complexes of $\Lambda$-modules modulo homotopical 
equivalence. Let $\D^b_{coh}(\Lambda)$ be the full subcategory of the 
bounded derived category of complexes of $\Lambda$-modules, consisting
of complexes whose cohomologies are of finite type. 
If $\Lambda$ is regular then the natural functor 
$\D_{parf}(\Lambda)\to \D^b_{coh}(\Lambda)$ is an equivalence 
of categories.~\QED
\end{Lm}

\smallskip
\noindent
1.3 Let $\cO\subset \EE$ be the ring of integers of $\EE$, $\kappa$ be
the residue field of $\cO$. Let $\bar R$ be a local noetherian
$\cO$-algebra with maximal ideal $\bar m$ and residue field $\kappa$. Assume
that $\bar R$ is complete in the $\bar \gm$-adic topology. Let
$\Sh(Y, \bar R)$ denote the category of $\bar m$-adic constructable
$\bar R$-sheaves (SGA5, V, 3.1.1). 
We say that $F\in \Sh(Y, \bar R)$ is \select{smooth of rank $m$} if for each $k>0$ 
the sheaf $F\otimes_{\bar R}\bar R/\bar \gm^k$ is locally constant 
and its fibres are free $\bar R/\bar \gm^k$-modules of rank $m$.

 Assume, in addition, that $\bar R$ is regular. Denote by $\D^b_c(Y,\bar R)$
the triangulated category defined in (\cite{E}, Theorem 6.3) (and denoted
$\D^b_c(Y-\bar R)$ in \select{loc. cit.}). So, $\D^b_c(Y,\bar R)$ 
is $\bar R$-additive,
has the natural $t$-structure, whose core
is equivalent to $\Sh(Y,\bar R)$, and the perverse $t$-structure
$p_{1/2}$ obtained by the gluing procedure (\cite{BBD}, 2.2.14). 
So, the usual definition of the perverse $t$-structure $p_{1/2}$ 
is applicable:

\begin{Def} An object $K\in\D^b_c(Y,\bar R)$ lies in $^p\D^{\leqslant 0}_c(Y,\bar R)$
(resp., in $^p\D^{\geqslant 0}_c(Y,\bar R)$) if and only if for any irreducible 
closed subscheme $Y'\subset Y$ there is a nonempty
open subscheme $U\subset Y'$ such that $\H^j(i^*_U K)=0$ for $j>-\dim U$ 
(resp., $\H^j(i_U^! K)=0$ for $j<-\dim U$), where $i_U:U\hook{}Y$ is the inclusion.
\end{Def}

 We let $\Perv(Y,\bar R)={^p\D^{\leqslant 0}_c(Y,\bar R)}\cap \,
{^p\D^{\geqslant 0}_c(Y,\bar R)}$. Any object of $\D^b_c(Y,\bar R)$ 
has a finite amplitude with respect to both the natural and the 
perverse $t$-structure $p_{1/2}$.
 
 The perverse $t$-structure $p_{1/2}$ on
$\D^b_c(Y,\bar R)$ is not preserved by Verdier duality (unless $\bar R=\kappa$).
This is already observed for $Y=\Spec k$. We have an 
equivalence of triangulated categories that preserves natural $t$-structures
as well as tensor products and internal $\Hom$'s (\cite{E}, 7.2)
$$
\D^b_c(\Spec k,\bar R)\iso \D^b_{coh}(\bar R),
$$  
and $\DD(\D^{\leqslant 0}_{coh}(\bar R))\subset \D^{\geqslant
0}_{coh}(\bar R)$. However, the image 
of $\D^{\geqslant 0}_{coh}(\bar R)$ under the Verdier duality functor $\DD$
is strictly bigger than $\D^{\leqslant 0}_{coh}(\bar R)$ 
(this image is described by \cite{GL}, 6.3.2).

\begin{Lm}
\label{Lm_Verdier_forever}
If $K\in{^p\D^{\leqslant 0}_c(Y,\bar R)}$ then 
$\DD K\in{^p\D^{\geqslant 0}_c(Y,\bar R)}$. 
\end{Lm}
\begin{Prf}
Let $Y'\subset Y$ be an irreducible closed subscheme. Chose a nonempty open subscheme
$U\subset Y'$ such that $i_U^*K$ is placed in usual degrees $\le -\dim U$ and all 
$\H^j(i_U^*K)$ are locally constant (by this we mean that 
$\H^j(i_U^*K)\otimes_{\bar R}\bar R/\bar\gm^k$ are locally constant for 
all $k>0$). Then $\DD(i_U^*K)$ is placed in usual degrees $\ge -\dim U$.
\end{Prf}

\medskip

We have a conservative triangulated functor of extension of scalars (\cite{E}, 6.3)
$$
\kappa\Lotimes_{\bar R}: \D^b_c(Y,\bar R)\to\D^b_c(Y,\kappa)
$$
  
\begin{Lm}
\label{Lm_affine_maps_bar_R}
 1) For $K\in \D^b_c(Y,\bar R)$ we have 
$K\in {^p\D^{\leqslant 0}_c(Y,\bar R)}$ if and only
if $\kappa\Lotimes_{\bar R}K\in {^p\D^{\leqslant 0}_c(Y,\kappa)}$.
Besides, if $\kappa\Lotimes_{\bar R}K\in {^p\D^{\geqslant 0}_c(Y,\kappa)}$
then $K\in {^p\D^{\geqslant 0}_c(Y,\bar R)}$.

\smallskip
\noindent
2) Let $f:Y\to Y'$ be a morphism of $k$-schemes of finite type. If $f$ is affine then
$f_*$ is right exact with respect to $p_{1/2}$. If $f$ is an open immersion then
$f_*$ (resp., $f_!$) is left exact (resp., right exact) with respect to $p_{1/2}$
{\rm (}the functors are from $\D^b_c(Y,\bar R)$ to $\D^b_c(Y',\bar R)${\rm ).} 
\end{Lm}
\begin{Prf}
1) By (\cite{E}, 6.3), the functors $i_U^*, i_U^!$ commute with 
$\kappa\Lotimes_{\bar R}$. To conclude, apply (\cite{E}, 3.1). \\
2) By (\select{loc.cit.}, 6.3), $f_*, f_!$ commute with $\kappa\Lotimes_{\bar R}$. 
It remains to use 1) and (\cite{BBD}, 4.1.1). 
\end{Prf}

\begin{Def} We say that $K\in\D^b_c(Y,\bar R)$ is
a \select{$\bar R$-flat perverse $\bar R$-sheaf} on $Y$ 
if it satisfies one of the following equivalent conditions:
\begin{itemize}
\item[i)] $K\in{^p\D^{\leqslant 0}_c(Y,\bar R)}$ and 
$\DD K\in{^p\D^{\leqslant 0}_c(Y,\bar R)}$;
\item[ii)] $K\in \Perv(Y,\bar R)$ and $\DD K\in\Perv(Y,\bar R)$;
\item[iii)] $\kappa\Lotimes_{\bar R}K\in \Perv(Y,\kappa)$
\end{itemize}
We denote by $\Perv_{fl}(Y,\bar R)$ the full subcategory of $\Perv(Y,\bar R)$
consisting of $\bar R$-flat perverse $\bar R$-sheaves on $Y$.
\end{Def}

 Note that the Verdier duality functor preserves
$\bar R$-flat perverse $\bar R$-sheaves and induces an autoequivalence 
of $\Perv_{fl}(Y,\bar R)$. 

\medskip
\noindent
1.4 \  Given a local homomorphism of $\cO$-algebras $\sigma: \bar R\to\cO$, denote
by $\gm_{\sigma}$ the kernel of the induced map $\bar
R\otimes_{\cO}\EE\to\EE$, and by $\bar R_{\sigma}$ the localization of
$\bar R\otimes_{\cO}\EE$ in $\gm_{\sigma}$. By Lemma~\ref{Lm_for_Gabber},
we have a triangulated category $\D^b_c(Y,\bar R)\otimes_{\bar R}\bar
R_{\sigma}$, which will also be denoted $\D^b_c(Y,\bar R)_{\sigma}$. It is
equiped with the natural and perverse $t$-structures induced by those of
$\D^b_c(Y,\bar R)$. So, by definition, 
$^p\D^{\leqslant 0}_c(Y,\bar R)_{\sigma}$ (resp., 
$^p\D^{\geqslant 0}_c(Y,\bar R)_{\sigma}$) is the essential image of
$^p\D^{\leqslant 0}_c(Y,\bar R)$ (resp., of $^p\D^{\geqslant 0}_c(Y,\bar R)$). 

\smallskip

 The core of the perverse $t$-structure on $\D^b_c(Y,\bar R)_{\sigma}$ 
will be denoted $\Perv(Y,\bar R)_{\sigma}$. 
The category $\Perv(Y,\cO)_{\sigma}$ will also be denoted by $\Perv(Y,\EE)$.
We also write $\Sh(Y,\EE)$ for $\Sh(Y,\cO)\otimes_{\cO}\EE$. 

\smallskip

 By 1) of Lemma~\ref{Lm_for_Gabber}, the functors 
defined for $\D^b_c(Y,\bar R)$ in \cite{E} extend trivially
to $\D^b_c(Y,\bar R)_{\sigma}$. So, for a
separated morphism $f:Y\to Y'$ of schemes of finite type over $k$ 
we have functors
$$
f_!, f_*: \D^b_c(Y,\bar R)_{\sigma}\to
\D^b_c(Y',\bar R)_{\sigma}
$$
and 
$$
f^!, f^*: \D^b_c(Y',\bar R)_{\sigma}\to
\D^b_c(Y, \bar R)_{\sigma}
$$
functorial in $f$. We also have triangulated functors
$$
(-)\Lotimes_{\bar R_{\sigma}}(-): \D^b_c(Y,\bar R)_{\sigma}\times
\D^b_c(Y,\bar R)_{\sigma} \to \D^b_c(Y,\bar R)_{\sigma}
$$
and 
$$
\R\HOM(-,-): \D^b_c(Y,\bar R)_{\sigma}^{op}\times \D^b_c(Y,\bar R)_{\sigma}
\to \D^b_c(Y,\bar R)_{\sigma}
$$
  
 The next lemma claims that the perverse $t$-structure on 
$\D^b_c(Y,\bar R)_{\sigma}$ defined by Lemma~\ref{Lm_for_Gabber} is, in fact,
obtained by a gluing procedure as in (\cite{BBD}, 2.2).
 
\begin{Lm} $K\in{^p\D^{\leqslant 0}_c(Y,\bar R)_{\sigma}}$ 
(resp., $K\in{^p\D^{\geqslant 0}_c(Y,\bar R)_{\sigma}}$)
if and only if
for any irreducible closed subscheme $Y'\subset Y$ there is a 
nonempty open subscheme $U\subset Y'$ such that $\H^j(i^*_U K)=0$ for
$j> -\dim U$ (resp., $\H^j(i^!_U K)=0$ for $j< -\dim U$), 
where $i_U:U\hook{} Y$ is the inclusion, and
the cohomologies are calculated with respect to the 
usual $t$-structure on $\D^b_c(U,\bar R)_{\sigma}$.
\end{Lm}
\begin{Prf}
The only if part is trivial. The if part follows from (\cite{BBD}, 1.4.12).
\end{Prf} 

\medskip

 From Lemma~\ref{Lm_Verdier_forever} and  
2) of Lemma~\ref{Lm_affine_maps_bar_R} we immediately get the next corolary.

\begin{Cor} 
\label{Cor_budu_primenjat'}
1) If $K\in {^p\D^{\leqslant 0}_c(Y,\bar R)_{\sigma}}$ then
$\DD K\in {^p\D^{\geqslant 0}_c(Y,\bar R)_{\sigma}}$.

\smallskip\noindent
2) Let $f:Y\to Y'$ be a morphism of $k$-schemes of finite type. If $f$ is affine then
$f_*$ is right exact with respect to $p_{1/2}$. If $f$ is an open immersion then
$f_*$ (resp., $f_!$) is left exact (resp., right exact) with respect to $p_{1/2}$
{\rm (}the functors are from $\D^b_c(Y,\bar R)_{\sigma}$ to 
$\D^b_c(Y',\bar R)_{\sigma}${\rm ).}  \QED
\end{Cor}

 Let $\bar R\to\bar R'$ be a local homomorphism of complete local
noetherian regular $\cO$-algebras with residue fields $\kappa$. It induces a conservative
triangulated functor of extension of scalars (\cite{GL}, A.1)
\begin{equation}
\label{ext_scalars_initial}
\bar R'\Lotimes_{\bar R}\; : \D^b_c(Y,\bar R)\to 
\D^b_c(Y,\bar R')
\end{equation} 
By (\select{loc.cit.}, A.1.5), the functors $f_!,f_*,f^*$ commute with 
(\ref{ext_scalars_initial})\footnote{The author does not know if 
$f^!$ and $\DD$ commute with (\ref{ext_scalars_initial}).}.

 Assume that $\sigma$ is the restriction to $\bar R$ of a
local homomorphism of $\cO$-algebras $\sigma': \bar R'\to\cO$.
Then $\bar R'\Lotimes_{\bar R}$ factors to give a triangulated 
functor 
\begin{equation}
\label{ext_scalars}
\D^b_c(Y,\bar R)_{\sigma}\to \D^b_c(Y,\bar R')_{\sigma'}
\end{equation}
that we denote by $\bar R'_{\sigma'}\Lotimes_{\bar R_{\sigma}}$. 
The map $\sigma$ itself yields the functor of 
extension of scalars 
$$
\EE\Lotimes_{{\bar R}_{\sigma}}: \D^b_c(Y,\bar R)_{\sigma}\to \D^b_c(Y, \EE),
$$ 
where, by definition, $\D^b_c(Y,\EE)$ is the triangulated category
$\D^b_c(Y,\cO)\otimes_{\cO}\EE$ (\cite{BBD}, 2.2.18).  

\begin{Lm} i) The functor $\EE\Lotimes_{\bar
R_{\sigma}}: \D^b_c(Y,\bar R)_{\sigma}\to
\D^b_c(Y, \EE)$ is triangulated and conservative.\\
ii) If for an object
$M$ of $\D^b_c(Y,\bar R)_{\sigma}$ we have
$\H^r(\EE\Lotimes_{\bar R_{\sigma}} M)=0$ for some $r$, then $\H^r(M)=0$,
where the cohomologies are calculated with respect to the usual
$t$-structures. 
\end{Lm}

\begin{Prf}
i) We must show that for an object $M$ of $\D^b_c(Y,\bar
R)_{\sigma}$ the property $\EE\Lotimes_{\bar
R_{\sigma}} M=0$ implies $M=0$. This follows from ii). \\
ii) Since for any closed point $i:\Spec k\to Y$ the functors
$\cO\Lotimes_{\bar R}$ and
$i^*$ commute, the assertion is fibrewise on $Y$, and we may assume
$Y=\Spec k$. Then we have an equivalence
$\D^b_c(Y,\bar R)\iso \D^b_{coh}(\bar R)\;$ (\cite{E}, 7.2). 
We have a commutative diagram of functors
$$
\begin{array}{ccc}
\D^b_{coh}(\bar R)\otimes_{\bar R}\bar R_{\sigma} & 
\to & \D^b_{coh}(\cO)\otimes_{\cO}\EE\\
\downarrow && \downarrow\\
\D^b_{coh}(\bar R_{\sigma}) & \to &
\D^b_{coh}(\EE),
\end{array}
$$
where the low horizontal arrow is the functor $\EE\Lotimes_{\bar R_{\sigma}}$
in the literal sense. By our assumption, the image of 
$\EE\Lotimes_{\bar R_{\sigma}} M$ under 
$\D^b_{coh}(\cO)\otimes_{\cO}\EE\to \D^b_{coh}(\EE)$ has no 
cohomology in degree $r$.

\begin{Slm}
\label{Slm_1}
Let $K$ be an object of $\D^b_{coh}(\bar R_{\sigma})$ such that 
$\H^r(\EE\Lotimes_{\bar R_{\sigma}}\! K)=0$ for some $r$ then 
$\H^r(K)=0$ (here $\EE\Lotimes_{\bar R_{\sigma}}\! K$ 
is understood in the literal sense).
\end{Slm}
\begin{Prf} Represent $K$ by a bounded complex $K_0$ of free $\bar
R_{\sigma}$-modules of finite type. Then $K_0=K_1\oplus K_2$, where $K_1$
is acyclic and the differential in $K_2$ is zero modulo the maximal 
ideal of $\bar R_{\sigma}$. Our assertion follows now from lemma of
Nakayama. 
\end{Prf}

\medskip

 By Sublemma~\ref{Slm_1}, the image of $M$ in $\D^b_{coh}(\bar R_{\sigma})$
has no cohomology in degree $r$. Represent $M$ by an object $M_0$ of
$\D^b_{coh}(\bar R)$. Since $\bar R_{\sigma}$ is flat over $\bar R$,
we have $\H^r(M_0\Lotimes_{\bar R}\bar R_{\sigma})\iso
\H^r(M_0)\otimes_{\bar R}\bar R_{\sigma}$ in the category of $\bar
R_{\sigma}$-modules. So, $\H^r(M_0)\otimes_{\bar R}\bar R_{\sigma}=0$
as $\bar R_{\sigma}$-module. Since $\H^r(M_0)$ is a
$\bar R$-module of finite type, the image of $\H^r(M_0)$ in 
$({\bar R}-mod)\otimes_{\bar R}\bar R_{\sigma}$ vanishes by
Remark~\ref{Rem_strange}. 
\end{Prf}

\medskip

 The perverse
$t$-structure on $\D^b_c(Y,\bar R)_{\sigma}$ is not preserved by Verdier duality
(unless $\bar R$ is of dimension 1). As in Sect. 1.3, we give the next definition.
\begin{Def} Let $\Perv_{fl}(Y,\bar R)_{\sigma}$ be the full subcategory of
$\D^b_c(Y,\bar R)_{\sigma}$ consisting of objects satisfying one of the following
equivalent conditions:
\begin{itemize}
\item[i)] $K\in{^p\D^{\leqslant 0}_c(Y,\bar R)_{\sigma}}$ and 
$\DD K\in{^p\D^{\leqslant 0}_c(Y,\bar R)_{\sigma}}$ ;
\item[ii)] $K\in \Perv(Y,\bar R)_{\sigma}$ and $\DD K\in\Perv(Y,\bar R)_{\sigma}$
\end{itemize}
The Verdier duality functor induces an autoequivalence of 
$\Perv_{fl}(Y,\bar R)_{\sigma}$. 
\end{Def} 

 Define the category of \select{smooth $\bar R_{\sigma}$-sheaves of rank $m$ on
$Y$} as the full subcategory of $\D^b_c(Y,\bar R)_{\sigma}$ consisting of $K$ 
such that $\EE\Lotimes_{\bar R_{\sigma}}K$ is a smooth $\EE$-sheaf of rank $m$ on $Y$ 
placed in usual cohomological degree zero. The next lemma will be used 
in the proof of Proposition~\ref{Pp_majus'}. 
 
\begin{Lm}
\label{Lm_for_Pp_majus'}
 1) Assume that $Y$ is connexe, and $y\to Y$ is a geometric point. 
Then the functor that sends $K$ to $K_y$ is an equivalence between the category
of smooth $\bar R_{\sigma}$-sheaves of rank $m$ on $Y$ and
the category of paires $(M,\rho)$, where $M$ is a free 
$\bar R_{\sigma}$-module of rank $m$, $\rho:\pi_1(Y,y)\to \Aut
M$ is a homomorphism such that there exists a finite type $\bar R$-submodule 
$M'\subset M$ generating $M$, invariant under $\pi_1(Y,y)$, and such that
the homomorphism $\rho:\pi_1(Y,y)\to \Aut M'$ is continuous (note that we do not
require $M'$ to be free over $\bar R$). \\
2) Let $U\subset Y$ be an open subscheme
and $K\in\D^b_c(Y,\bar R)_{\sigma}$. If both $K\Lotimes_{\bar R_{\sigma}}\EE$ and
$(\DD K)\Lotimes_{\bar R_{\sigma}}\EE$ are perverse and the Goresky-MacPherson
extensions from $U$ to $Y$ then $K\in\Perv_{fl}(Y,\bar R)_{\sigma}$ 
is the Goresky-MacPherson extension from $U$ to $Y$.
\end{Lm}
\begin{Prf} 1) Argue as in SGA5, VI, 1.2.\\
2) For $M\in \D^b_c(Y,\bar R)_{\sigma}$
the condition $M\Lotimes_{\bar R_{\sigma}}\EE\in{^p\D^{\leqslant 0}_c(Y,\EE)}$
implies $M\in {^p\D^{\leqslant 0}_c(Y,\bar R)_{\sigma}}$. So, our assumption
implies $K\in\Perv_{fl}(Y,\bar R)_{\sigma}$. Let $i:S\to Y$ be a closed subscheme
whose complement is $U$. The complexes $i^*(K\Lotimes_{\bar R_{\sigma}}\EE)$
and $i^*((\DD K)\Lotimes_{\bar R_{\sigma}}\EE)$ lie in $^p\D^{< 0}_c(S,\EE)$. 
So, $i^*K$ and $\DD i^!K$ lie in $^p\D^{< 0}_c(S,\bar R)_{\sigma}$. By
Corolary~\ref{Cor_budu_primenjat'}, $i^!K\in{^p\D^{>0}_c(S,\bar R)_{\sigma}}$.
\end{Prf}

\medskip\smallskip
\noindent
1.5 \select{Laumon's sheaf $\cL^d_{\bar E}$} 

\smallskip
\noindent
Let $\bar R$ be as in Sect. 1.3. Assume that $\bar R$ is of
characteristic zero. 
Let $\bar E$ be a smooth $\bar R$-sheaf 
on $X$. Denote by $\sym: X^d\to X^{(d)}$ the natural map. 
Consider the smooth $\bar R$-sheaf
$\bar E^{\boxtimes\, d}$ on $X^d$ (the tensoring is taken over $\bar R$).
Set
$$
\bar E^{(d)}=(\sym_!(\bar E^{\boxtimes\, d}))^{S_d}
$$ 
This is a direct summand of the constructable $\bar R$-sheaf $\sym_!(\bar
E^{\boxtimes\, d})$ on $X^{(d)}$. Since $\sym_!(\bar E^{\boxtimes\, d})[d]$
is a $\bar R$-flat perverse $\bar R$-sheaf, which is the Goresky-MacPherson 
extension of its restriction to any nonempty open subscheme, the same holds
for $\bar E^{(d)}[d]$.   

  Following \cite{L1}, we associate to $\bar E$ a perverse 
$\bar R$-sheaf $\cL^d_{\bar E}$ on $\Sh^d_0$ that we call Laumon's sheaf.   
Denote by $\Fl^{1,\ldots,1}$ (1 occurs $d$ times) the stack of complete 
flags $(F_1\subset\ldots\subset F_d)$, where $F_i$ is a coherent torsion 
sheaf on $X$ of length $i$. The morphism $\gp_0: \Fl^{1,\ldots,1}\to\Sh^d_0$ 
that sends $(F_1\subset\ldots\subset F_d)$ to $F_d$ is representable and 
proper. The morphism
$\gq_0:\Fl^{1,\ldots1}\to\Sh^1_0\times\ldots\times\Sh^1_0$ that sends
$(F_1\subset\ldots\subset F_d)$ to $(F_1, F_2/F_1,\ldots, F_d/F_{d-1})$ 
is a generalized affine fibration. 
This, in particular, implies that $\Fl^{1,\ldots,1}$ is smooth. 

 Springer's sheaf $\Spr^d_{\bar E}$ is defined as
$$
\Spr^d_{\bar E}=(\gp_0)! (\gq_0)^* (\div^{\times d})^* (\bar E^{\boxtimes\, d})
$$
Laumon's theorem claims that $\gp_0$ is \select{small}. 
It follows that $\Spr^d_{\bar E}$ is a $\bar R$-flat perverse $\bar R$-sheaf, 
which is the Goresky-MacPherson extension of its restriction to any nonempty 
open substack. It also carries a natural action of the symmetric group $S_d$ 
(\cite{L1}, 3.3.1). Set
$$
\cL^d_{\bar E}=(\Spr^d_{\bar E})^{S_d},
$$
where the invariants are taken in $\Perv(\Sh^d_0,\bar R)$. 
As a direct summand of $\Spr^d_{\bar E}$, the perverse $\bar R$-sheaf 
$\cL^d_{\bar E}$ is $\bar R$-flat and coincides with the Goresky-MacPherson extension
of its restriction to any nonempty open substack. Besides, the formation of 
$\cL^d_{\bar E}$ commutes with extension of scalars 
(\ref{ext_scalars_initial}), and the Verdier dual 
of $\cL^d_{\bar E}$ is canonically isomorphic to $\cL^d_{\bar E^*}$.

\section{Automorphic sheaves $\Aut_{\bar E}$}

2.1 Let $\bar R,\bar E$ be as in Sect. 1.5. 
Assume that $\EE$ contains the group of $p$-th roots of unity, 
so that we can fix a nontrivial additive character $\psi:\Fp\to \cO^*$. 
Then the Artin-Schreier 
sheaf $\cL_{\psi}$ associated to $\psi$ can be viewed as a smooth $\cO$-sheaf 
(or smooth $\bar R$-sheaf) of rank 1 on $\A^1$. 

 Let $n>0,d\ge 0$. In Sect 2.1 and 4.1 of \cite{Ly2} 
we introduced the diagram
$$
\begin{array}{ccccc}
^{\le n}{\Sh^d_0} & \getsup{\beta} & {_n\cQ_d} & \toup{\mu}& \A^1\\
&& \downarrow\lefteqn{\scriptstyle \varphi} & \searrow\lefteqn{\scriptstyle \zeta}\\
&&  _n\cY_d & \to & _n\cM_d
\end{array}
$$  
\begin{Def} We associate to $\bar E$ an object 
$_n\cF^d_{\bar E,\psi}\in \Perv_{fl}({_n\cQ_d}, \bar R)$ given by
$$
_n\cF^d_{\bar E,\psi}=\beta^*\cL^d_{\bar E}\otimes_{\bar R} \mu^*\cL_{\psi}
[b](\frac{b}{2}),
$$
where $b=\dim{_n\cQ_d}$. We also define 
$_n\cP^d_{\bar E,\psi}\in \D^b_c({_n\cY_d},\bar R)$ and 
$_n\cK^d_{\bar E}\in \D^b_c({_n\cM_d}, \bar R)$ by
$_n\cP^d_{\bar E,\psi}=\varphi_!(_n\cF^d_{\bar E,\psi})$ and
$\;{_n\cK^d_{\bar E}}=\zeta_!(_n\cF^d_{\bar E,\psi})$.
\end{Def}

 The formation of all these complexes
commutes with extension of scalars (\ref{ext_scalars_initial}).  
The complex $_n\cK^d_{\bar E}$
does not depend on $\psi$ in the following sense.

\begin{Lm} For any two nontrivial additive characters $\psi,\psi':\Fp\to \cO^*$
there is a canonical isomorphism $\zeta_!(_n\cF^d_{\bar E,\psi})\iso
\zeta_!(_n\cF^d_{\bar E,{\psi'}})$ in $\D^b_c({_n\cM_d},\bar R)$. 
\end{Lm}
\begin{Prf}
There is a unique $a\in\Fp^*$ such that $\psi'(x)=\psi(ax)$ for $x\in\Fp$. So, 
$\cL_{\psi'}\iso a^*\cL_{\psi}$, where $a:\A^1\iso\A^1$ denotes the multiplication by
$a$. Let $\alpha$ denote the automorphism of $_n\cQ_d$ that multiplies 
$s_i:\Omega^{n-i}\iso L_i/L_{i-1}$ by $a^{i-1}$ for $i=2,\ldots,n$, where 
$(L_1\subset\ldots\subset L_n\subset L, (s_i))$ is a point of $_n\cQ_d$. Then
$\alpha^*(_n\cF^d_{\bar E,\psi})\iso {_n\cF^d_{\bar E,\psi'}}$, and our assertion
follows.
\end{Prf}

\medskip
\noindent
2.2  \  Fix $\sigma:\bar R\to \cO$ as in Sect. 1.4.
Let $E_0$ be the image of $\bar E\otimes_{\bar R}\cO$ in $\Sh(X,\EE)$.
It was shown in (\cite{FGV2}, 7.3 and 7.5) that $_n\cP^d_{E_0,\psi}$ (resp.,
$_n\cK^d_{E_0}$) satisfies Hecke property
with respect to $E_0$. Let us show that this Hecke 
property still holds in the corresponding categories 
$\D^b_c(\; .\; ,\bar R)_{\sigma}$. 

\smallskip

 As for $_n\cP^d_{\bar E,\psi}$, 
(Proposition 5, Corolary 2 and Lemma 12, \cite{Ly2}) hold with
$E_0$ replaced by $\bar E$, if the maps and isomorphisms are understood as such
in $\D^b_c(\; .\; ,\bar R)_{\sigma}$. Indeed, one constructs the morphisms
in the same way and checks that they are isomorphisms applying the conservative
functor $\EE\Lotimes_{\bar R_{\sigma}}$. 

\smallskip

 The Hecke property of $_n\cK^d_{\bar E}$ is formulated as follows.
Let $_n\Mod_d$ denote the stack
classifying modifications $(L\subset L')$ of rank $n$ vector bundles on $X$
with $\deg(L'/L)=d$. Let $\supp:{_n\Mod_d}\to X^{(d)}$ be the map that sends
$(L\subset L')$ to $\div(L'/L)$. For $d'\ge 0$ denote by 
$$
\gp_{\cM}:{_n\cM_d}\times_{\Bun_n}{_n\Mod_{d'}}\to {_n\cM_{d+d'}}
$$ 
the map that sends
$(\Omega^{n-1}\hook{}L\hook{} L')$ to the composition $(\Omega^{n-1}\hook{} L')$. 
It is representable and proper. Let also 
$\gq_{\cM}:{_n\cM_d}\times_{\Bun_n}{_n\Mod_{d'}}\to{_n\cM_d}$ denote the projection.

\begin{Pp} 
\label{Pp_Hecke_for_cK}
For any smooth $\bar R$-sheaf $\bar E$ on $X$ and any $d\ge 0$ there is
a natural morphism
\begin{equation}
\label{iso_nuzhnyj}
(\gq_{\cM}\times\supp)_!\,\gp_{\cM}^*(_n\cK^{d+1}_{\bar E})\to {_n\cK^d_{\bar E}}\boxtimes
\bar E(\frac{2-n}{2})[2-n]
\end{equation}
in $\D^b_c({_n\cM_d}\times X, \bar R)_{\sigma}$, 
which is an isomorphism if $\rk \bar E\le n$.
\end{Pp}

 Let $_n\wt{\Mod}_d$ be the stack of flags $(L_0\subset\ldots\subset L_d)$, where
$(L_i\subset L_{i+1})\in{_n\Mod_1}$ for all $i$. Let $\wt{\supp}:{_n\wt{\Mod}_d}\to X^d$
be the map that sends $(L_0\subset\ldots\subset L_d)$ to
$(\div(L_1/L_0),\ldots,\div(L_d/L_{d-1}))$. Let 
$\tilde\gp_{\cM}:{_n\cM_d}\times_{\Bun_n}{_n\wt{\Mod}_{d'}}\to {_n\cM_{d+d'}}$
denote the composition 
$$
{_n\cM_d}\times_{\Bun_n}{_n\wt{\Mod}_{d'}}\to {_n\cM_d}\times_{\Bun_n}{_n\Mod_{d'}}
\toup{\gp_{\cM}} {_n\cM_{d+d'}},
$$
where the first map is the projection. We also have the map
$\gq_{\cM}\times\wt{\supp}: {_n\cM_d}\times_{\Bun_n}{_n\wt{\Mod}_{d'}}\to {_n\cM_d}\times
X^{d'}$.

\begin{Cor}
For any smooth $\bar R$-sheaf $\bar E$ on $X$ and any $d,d'\ge 0$ there is
a natural morphism 
\begin{equation}
\label{arrow_hecke_1}
(\gq_{\cM}\times\wt{\supp})_!\,\tilde\gp_{\cM}^*(_n\cK^{d+d'}_{\bar E})\to 
{_n\cK^d_{\bar E}}\boxtimes \bar E^{\boxtimes d'}(\frac{2d'-nd'}{2})[2d'-nd']
\end{equation}
in $\D^b_c({_n\cM_d}\times X^{d'}, \bar R)_{\sigma}$, 
which is an isomorphism if $\rk\bar E\le n$. \QED
\end{Cor}

As in (\cite{Ly2}, Sect 6.6) we write $^{rss}X^{d'}$   
for the open subscheme of $X^{d'}$ that parametrizes pairwise different points
$(x_1,\ldots,x_{d'})\in X^{d'}$ (`rss' stands for `regular semisimple'). We also denote by
$_n^{rss}\wt{\Mod}_{d'}$ the preimage of $^{rss}X^{d'}$ under $\wt{\supp}$. The symmetric
group $S_{d'}$ acts on 
$$
_n\cM_d\times_{\Bun_n}{_n^{rss}\wt{\Mod}_{d'}},
$$ 
and this action lifts naturally 
to an action on $\tilde\gp_{\cM}^*(_n\cK^{d+d'}_{\bar E})$. Since the restriction 
${_n^{rss}\wt{\Mod}_{d'}}\to {^{rss}X^{d'}}$ of $\wt{\supp}$ is $S_{d'}$-equivariant,
$S_{d'}$ acts on the complex 
$$
(\gq_{\cM}\times\wt{\supp})_!\,\tilde\gp_{\cM}^*(_n\cK^{d+d'}_{\bar E})
$$
restricted to $_n\cM_d\times {^{rss}X^{d'}}$. On the other hand, $S_{d'}$ acts on
${\bar E}^{\boxtimes d'}$ and, hence, on the right hand side of (\ref{arrow_hecke_1}).
Using the explicit description of the map (\ref{arrow_hecke_1}) (cf. Sect. 2.3)
one easily obtains the next result.

\begin{Lm} 
\label{Lm_Hecke_equivariance}
The map (\ref{arrow_hecke_1})
restricted to $_n\cM_d\times {^{rss}X^{d'}}$ is $S_{d'}$-equivariant. \QED
\end{Lm}

\noindent
2.3 \  In this subsection we prove Proposition~\ref{Pp_Hecke_for_cK}.

 Denote by $_n\tilde\cY^{1}_d$ the stack of collections $((t_i),L,
x\in X)$, where $L$ is a rank $n$ vector bundle on $X$ with $\deg
L=d+\deg(\Omega^{(n-1)+\ldots+(n-n)})$ and 
$$
t_i:\Omega^{(n-1)+\ldots+(n-i)}\hook{} 
(\wedge^i L)(x)
$$ 
are sections satisfying Pl\"ucker's relations as in (\cite{Ly2}, Sect. 4.1).  
We have a closed immersion $_n\cY_d\times X\hook{} {_n\tilde\cY^1_d}$ given by 
the condition: every $t_i$ factors as 
$\Omega^{(n-1)+\ldots+(n-i)}\hook{}\wedge^i L\subset (\wedge^i L)(x)$.

\smallskip

 If $(L\subset L')\in{_n\Mod_1}$ with $x=\div(L'/L)$ then $(\wedge^i L')(-x)\subset 
\wedge^i L\subset \wedge^i L'$ for all $i$. This allows to define a map
$$
\tilde\gq_{\cY}: {_n\cY_{d+1}}\times_{\Bun_n}{_n\Mod_1}\to {_n\tilde\cY^1_d}
$$ 
that sends
$(((t'_i), L')\in {_n\cY_{d+1}}, L\subset L')$ to $((t_i),L,x)$, where $x=\div(L'/L)$
and $t_i$ is the composition 
$$
t_i: \Omega^{(n-1)+\ldots+(n-i)}\hook{t'_i} \wedge^i L'\subset (\wedge^i L)(x)
$$
The map $\tilde\gq_{\cY}$ is representable and proper. 

\begin{Def} For any smooth $\bar R$-sheaf $\bar E$ on $X$ set
$_n\tilde\cP^{d,1}_{\bar E,\psi}=(\tilde\gq_{\cY})_!(_n\cP^{d+1}_{\bar E,\psi}\boxtimes
\bar R)(\frac{n}{2})[n]$.
\end{Def}

\begin{Rem}
i) It may be shown that $_n\tilde\cP^{d,1}_{\bar E,\psi}$ lies in 
$\Perv({_n\tilde\cY^{1}_d},\bar R)_{\sigma}$ and coincides with the 
Goresky-MacPherson extension of its restriction
to any nonempty open substack. Besides, the Verdier dual to 
$_n\tilde\cP^{d,1}_{\bar E,\psi}$ is canonically isomorphic to $_n\tilde\cP^{d,1}_{\bar
E^*,\psi^{-1}}$. We will not need these facts. \\
ii) The following square is cartesian
$$
\begin{array}{ccc}
_n\cY_d\times_{\Bun_n}{_n\Mod_1} & \hook{} & {_n\cY_{d+1}}\times_{\Bun_n}{_n\Mod_1}\\
\downarrow\lefteqn{\scriptstyle \gq_{\cY}\times\supp} &&
\downarrow\lefteqn{\scriptstyle \tilde\gq_{\cY}}\\
_n\cY_d\times X & \hook{} & _n{\tilde\cY}^1_d
\end{array}
$$
So, the restriction of 
$_n\tilde\cP^{d,1}_{\bar E,\psi}$ to $_n\cY_d\times X$ is described by by 
(\cite{Ly2}, Prop.~5).
\end{Rem}

 The only property of $_n\tilde\cP^{d,1}_{\bar E,\psi}$ we need is the following.
Let $_n'{\tilde\cY}^1_d\hook{}{_n{\tilde\cY}^1_d}$ be the closed substack given by:
$t_1$ factors as $\Omega^{n-1}\hook{}L\hook{}L(x)$. So, $_n\cY_d\times X\hook{}
{_n'{\tilde\cY}^1_d}$ is a closed substack. Proposition~\ref{Pp_Hecke_for_cK}
will follow from the next observation.

\begin{Lm} 
\label{Lm_supp_tilde_cP}
For any smooth $\bar R$-sheaf $\bar E$ on $X$ 
the restriction of $_n\tilde\cP^{d,1}_{\bar E,\psi}$ to 
$_n'{\tilde\cY}^1_d$ vanishes outside $_n\cY_d\times X$. 
\end{Lm}
\begin{Prf} In the case $\bar R=\cO$ this follows from 
(\cite{FGV2}, 7.5) (the assumption $\rk E=n$ required in \select{loc.cit.}
is, in fact, not used for this particular statement).
Applying the conservative functor $\EE\Lotimes_{\bar R_{\sigma}}$,
one reduces the general case to $\bar R=\cO$.
\end{Prf}

\medskip

\begin{Prf}\select{of Proposition~\ref{Pp_Hecke_for_cK}} \ 
Define the closed substack $\cT$ of $_n\cY_{d+1}\times_{\Bun_n}{_n\Mod_1}$  from the
cartesian square
$$
\begin{array}{ccc}
\cT & \hook{} & _n\cY_{d+1}\times_{\Bun_n}{_n\Mod_1}\\
\downarrow && \downarrow\lefteqn{\scriptstyle \tilde\gq_{\cY}}\\
_n'{\tilde\cY}^1_d & \hook{}& {_n{\tilde\cY}^1_d}
\end{array}
$$
Then we have a commutative diagram, where the right square is cartesian
$$
\begin{array}{cccccccc}
_n\cY_d\times X & \hook{} & _n'{\tilde\cY}^1_d & \gets & \cT & \to & _n\cY_{d+1}\\
 & \searrow & \downarrow\lefteqn{\scriptstyle \xi_{\cY}} && \downarrow &&
\downarrow\\
&& _n\cM_d\times X & \getsup{\gq_{\cM}\times\supp} &
_n\cM_d\times_{\Bun_n}{_n\Mod_1} & \toup{\gp_{\cM}} & _n\cM_{d+1}
\end{array}
$$
We have denoted here by $\xi_{\cY}$ the natural forgetful map. Denote for brevity by
$\cF$ the restriction of $_n\tilde\cP^{d,1}_{\bar E,\psi}$ to $_n'{\tilde\cY}^1_d$. 
By the base change theorem, 
$$
(\gq_{\cM}\times\supp)_!\,\gp_{\cM}^*(_n\cK^{d+1}_{\bar E})\iso
(\xi_{\cY})_!\cF[-n](\frac{-n}{2})
$$
By Lemma~\ref{Lm_supp_tilde_cP}, $\cF$ is the extension by zero from 
$_n\cY_d\times X$. By (\cite{Ly2}, Prop. 5), we get a morphism
$\cF\to {_n\cP^d_{\bar E,\psi}}\boxtimes \bar E(1)[2]$, which is an isomorphism
when $\rk \bar E\le n$. Our assertion follows.
\end{Prf}

\medskip\smallskip
\noindent
2.4 \  Recall the definition of a Hecke-eigensheaf (\cite{FGV2}, 1.1). 
Consider the correspondence
\begin{equation}
\label{diag_Hecke}
\Bun_n\times X\getsup{\;\;\,\gq\times\supp}\;{_n\Mod_1}\toup{\gp}\,\Bun_n,
\end{equation}
where $\gp$ sends $(L\subset L')\in{_n\Mod_1}$ to $L'$, 
$\gq$ sends $(L\subset L')$ to $L$, and $\supp:{_n\Mod_1}\to X$
is the map defined in Sect. 2.2. 

 The \select{Hecke functor}
$\H^1_n: \D^b_c(\Bun_n,\bar R)_{\sigma}\to \D^b_c(\Bun_n\times X,\bar R)_{\sigma}$ 
is defined by
\begin{equation}
\label{hecke_functor}
\H^1_n(K)=(\gq\times\supp)_!\gp^*(K)(\frac{n-1}{2})[n-1]
\end{equation}
Consider the $d$-th iteration of $\H^1_n$:
$$
(\H^1_n)^{\boxtimes d}:\D^b_c(\Bun_n,\bar
R)_{\sigma}\to \D^b_c(\Bun_n\times X^d, \bar R)_{\sigma}
$$ 
For any $K\in\D^b_c(\Bun_n, \bar R)_{\sigma}$ the restriction of
$(\H^1_n)^{\boxtimes d}(K)$ to $\Bun_n\times\, {^{rss}X^d}$ 
is naturally equivariant with respect to the action of the symmetric
group $S_d$ on $^{rss}X^d$. 

\begin{Def}
\label{Def_Hecke_eigensheaf}
 Let $\bar E$ be a smooth $\bar R$-sheaf of rank $n$ on $X$. 
A \select{Hecke eigensheaf} with respect to $\bar E$ is a nonzero object
$K\in\D^b_c(\Bun_n,\bar R)_{\sigma}$ equiped with an isomorphism 
$\H^1_n(K)\iso K\boxtimes \bar E$ such that the resulting map
$$
(H^1_n)^{\boxtimes d}(K)\mid_{\Bun_n\times\, {^{rss}X^d}}\iso K\boxtimes\bar
E^{\boxtimes d}\mid_{\Bun_n\times\, {^{rss}X^d}}
$$
is $S_d$-equivariant.
\end{Def}

Following \cite{FGV2}, pick 
a line bundle $\cL^{est}$ on $X$ such that for any vector bundle $M$ on $X$ of rank
$k<n$, $\Hom(M,\cL^{est})=0$ implies that
\begin{itemize}
\item[a)] $\deg(M)>nk(2g-2)$,
\item[b)] $\Ext^1(\Omega^{k-1},M)=0$.
\end{itemize}
For example, $\cL^{est}$ of degree $>2n(2g-2)$ satisfies this property. Let us denote
by $_n\cM_d^{est}\subset{_n\cM_d}$ the open substack consisting of
$(\Omega^{n-1}\hook{}L)\in {_n\cM_d}$ such that $\Hom(L/\Omega^{n-1},\cL^{est})=0$.
Set $_n\cM=\cup_{d\ge 0}\;{_n\cM_d}$. Denote by $\varrho:{_n\cM\to\Bun_n}$
the natural projection. 

\medskip
\noindent
{\bf Notational Convention.}
For notational convenience, in what follows by \select{degree}
of a coherent sheaf on $X$ of generic rank $n$ we will understand its usual degree
$-n(n-1)(g-1)$, so that $\cO\oplus\Omega\oplus\ldots\oplus\Omega^{n-1}$ is of degree
zero. We write $\Bun_n^d$ for the connected component of $\Bun_n$ classifying
vector bundles of rank $n$ and degree $d$ on $X$.

\medskip

 Recall that $\bar R$ is a local complete noetherian regular $\cO$-algebra
with residue field $\kappa$, and $\sigma:\bar R\to\cO$ is a local homomophism
of $\cO$-algebras (in particular, $\bar R$ is of characteristic zero).
The geometric Langlands conjecture (relative to $\bar R$) may be formulated as 
follows (cf. \cite{L1,FGKV, FGV2}). 

\begin{Con}
\label{Con_1}
Let $\bar E$ be a smooth $\bar R$-sheaf on $X$ of rank $n$ such that
$E_0\otimes_{\EE}\Qlb$ is irreducible. Then for any $d\ge 0$ the natural map 
$\zeta_!({_n\cF^d_{\bar E,\psi}})\to \zeta_*({_n\cF^d_{\bar E,\psi}})$
is an isomorphism over $_n\cM_d^{est}$ in $\D^b_c({_n\cM_d^{est}},\bar R)_{\sigma}$, 
and the restriction of $_n\cK^d_{\bar E}$ to
$_n\cM_d^{est}$ lies in $\Perv_{fl}({_n\cM_d^{est}},\bar R)_{\sigma}$. 
Moreover, there exist the following data:
\begin{itemize}
\item a perverse sheaf $\Aut_{\bar E}\in\Perv_{fl}(\Bun_n,\bar R)_{\sigma}$ 
equiped with the structure of a Hecke-eigensheaf with respect to $\bar E$; 
\item for each $d\ge 0$ an isomorphism in $\D^b_c({_n\cM_d},\bar
R)_{\sigma}$ between $_n\cK^d_{\bar E}$ and the inverse image  
$$
\varrho^*\Aut^d_{\bar E}[d-n^2(g-1)](\frac{d-n^2(g-1)}{2})
$$ 
under $\varrho: {_n\cM_d}\to \Bun^d_n$ (we write $\Aut^d_{\bar E}$ for the
restriction of $\Aut_{\bar E}$ to $\Bun_n^d$). 
\end{itemize}
These data satisfy the following properties:
\begin{itemize}
\item[i)] the Hecke properties of $\Aut_{\bar E}$ and of $_n\cK_{\bar E}$ 
are compatible.
\item[ii)] $\Aut^d_{\bar E}$ is the Goresky-MacPherson extension from any nonempty
open substack of $\Bun^d_n$. 

\item[iii)] $\Aut_{\bar E}$ is cuspidal in the sense that for any
nontrivial partition $\bar n=(n_1,\ldots,n_k)$ of $n$ determining
the correspondence
$$
\begin{array}{ccc}
\Flag_{\bar n} & \toup {\gp_{\bar n}} & \Bun_n\\
\downarrow\lefteqn{\scriptstyle \gq_{\bar n}}\\
\Bun_{n_1}\times\ldots\times\Bun_{n_k},
\end{array}
$$
we have $(\gq_{\bar n})_!(\gp_{\bar n})^*\Aut_{\bar E}=0$. Here
$\Flag_{\bar n}$ is the stack of flags $(M_1\subset\ldots\subset
M_k)$, where $M_i/M_{i-1}\in \Bun_{n_i}$. The map $\gp_{\bar n}$
sends this flag to $M_k$, and $\gq_{\bar n}$ sends this flag to
$(M_1, \, M_2/M_1,\ldots, M_k/M_{k-1})$. 
\end{itemize}
\end{Con} 

\begin{Rems}
1) The property i) means the following. For any $d\ge 0$ we have a commutative 
diagram, where the left square is cartesian
$$
\begin{array}{ccccc}
_n\cM_d\times X & \getsup{\gq_{\cM}\times\supp} & _n\cM_d\times_{\Bun_n}{_n\Mod_1} &
\toup{\gp_{\cM}} & _n\cM_{d+1}\\
\downarrow && \downarrow && \downarrow\\
\Bun_n\times X & \getsup{\gq\times\supp} & _n\Mod_1 & \toup{\gp} & \Bun_n
\end{array}
$$ 
It is required that the restriction of the isomorphism 
$\H^1_n(\Aut_{\bar E})\iso \Aut_{\bar E}\boxtimes\bar E$ under the left 
vertical arrow in the above diagram coincides with the isomorphism 
(\ref{iso_nuzhnyj}) (up to a cohomological shift and a Tate twist).

\smallskip
\noindent
2) From main result of \cite{FGV2} it follows that if certain vanishing conjecture 
(Conjecture~2.3, \select{loc.cit}) is true then Conjecture~\ref{Con_1} is 
true for $\bar R=\cO$. Indeed, the first assertion follows from 
(3.6, \select{loc.cit.}) combined with the properties of the Fourier transform 
(\cite{L2}, 1.3.1.1, 1.3.2.3).

\smallskip
\noindent
3) The perverse sheaf $\Aut_{\bar E}$ in Conjecture~\ref{Con_1}
is defined up to a canonical isomorphism.
\end{Rems}

\medskip\smallskip\noindent
2.5 \ \select{The quotient of $\Bun_n$ by a 2-action of $\Gm$}

\smallskip\noindent
Consider the $k$-prestack whose category fibre at a scheme $S$ is the
following groupo\"\i d. Its objects are vector bundles $L$ on $S\times X$ of
rank $n$. A morphism from $L_1$ to $L_2$ is an equivalence
class $\in
\{(\cA,u)\}/
\sim$, where $\cA$ is an invertible sheaf on $S$, $u: L_1\iso
L_2\otimes\cA$ is an isomorphism of $\cO_{S\times X}$-modules, and the 
pairs $(\cA,u)$ and $(\cA',u')$ are equivalent if there exists an
isomorphism $\cA\iso \cA'$ making commute the diagram
$$
\begin{array}{ccc}
L_1 & \toup{u} & L_2\otimes\cA\\
& \searrow\lefteqn{\,\scriptstyle{u'}} & \downarrow\lefteqn{\wr}\\
&& L_2\otimes\cA'
\end{array}
$$
We define $\Bunb_n$ as the stack associated to this
prestack.\footnote{Scalar automorphisms of vector bundles provide
a 2-action of $\Gm$ on $\Bun_n$, that is, an action by
2-automorphisms of stacks. The quotient of a 1-stack under a 2-action
is, in general, a 2-stack. In our particular case this 2-stack is
representable by the 1-stack $\Bunb_n$.} Then the
natural morphism $\gr:\Bun_n\to \Bunb_n$ is a $\Gm$-gerb.

\begin{Lm}
$\Bunb_n$ is an algebraic stack locally of finite type and
smooth of pure dimension $n^2(g-1)+1$.
\end{Lm}
\begin{Prf}
Let $S$ be a noetherian scheme, $L_1$ and $L_2$ be vector bundles on $S\times X$ 
of rank $n$.
Consider the morphism $\Isom(L_1,L_2)\to S$ obtained by the base change
$S\toup{(L_1,L_2)}\Bun_n\times\Bun_n$ from the diagonal mapping
$\Bun_n\to\Bun_n\times\Bun_n$. Recall that $\Isom(L_1,L_2)$ is an open subscheme 
of some affine $S$-scheme $\V(\cF)$, where $\cF$ is a coherent $\cO_S$-module, 
and $\Isom(L_1,L_2)$ is of finite type over $S$ (\cite{LMb}, the proof of 4.6.2.1). 
We have a free action of $\Gm$ on 
$\Isom(L_1,L_2)$, and the square is cartesian
$$
\begin{array}{ccc}
\Bun_n\times_{\Bunb_n}\Bun_n & \to & \Bun_n\times\Bun_n\\
\uparrow && \uparrow\lefteqn{(L_1,L_2)}\\
\Isom(L_1,L_2)/\Gm & \to & S
\end{array}
$$
It follows that $\Isom(L_1,L_2)/\Gm$ is an open subscheme of $\P(\cF)$, in
particular it is separated over $S$. So, the diagonal mapping $\Bunb_n\to
\Bunb_n\times\Bunb_n$ is representable, separated and quasi-compact. 
 
  If $Y\to\Bun_n$ is a presentation of $\Bun_n$ then the composition
$Y\to \Bun_n\to\Bunb_n$ is a presentation of $\Bunb_n$. 
\end{Prf}

\medskip

The connected components of $\Bunb_n$ are numbered by $d\in\ZZ$:
the component $\Bunb^d_n$ is the image of $\Bun^d_n$ under $\gr:\Bun_n\to \Bunb_n$. 
The morphism $\gr:\Bun^d_1\to \Bunb^d_1$
is, in fact, the canonical map $\Pic^d X\to\uPic^d X$, where
$\uPic^d X$ is the Picard scheme of $X$.

 Let $\mult:\Bun_n\times X\to\Bun_n$ be the map that sends $(L,x)$ to $L(x)$. 
There is a unique map 
$\multb:\Bunb_n\times X\to\Bunb_n$ making commute the diagram
\begin{equation}
\label{diag_mult}
\begin{array}{ccc}
\Bun_n\times X & \toup{\mult} & \Bun_n\\
\downarrow\lefteqn{\scriptstyle \gr\times\id} && 
\downarrow\lefteqn{\scriptstyle \gr}\\
\Bunb_n\times X & \toup{\multb} & \Bunb_n
\end{array}
\end{equation}
The composition $_1\cM_d\to\Bun_1^d\to\Bunb_1^d$ is, in fact, 
the Abel-Jacobi map $X^{(d)}\to \uPic^d X$ sending a divisor $D$ to the
isomorphism class of $\cO(D)$. More generally,  
the composition $_n\cM_d\to\Bun^d_n\toup{\gr}\Bunb^d_n$ is 
representable. Given a scheme $S$ and an $S$-point
$S\toup{L}\Bun^d_n$, consider the scheme $Z_S={_n\cM_d}\times_{\Bun_n}
S$. This is an $S$-scheme that classifies nonzero sections 
$s: \Omega^{n-1}\hook{} L$. The group $\Gm$ acts freely on $Z_S$ 
(over $S$), multiplying $s$ by a scalar. One checks that 
$_n\cM_d\times_{\Bunb_n} S$ is identified with the quotient $Z_S/\Gm$.  
In particular, for $S=\Spec k$ the scheme $Z_S$ is the projective space 
$(\Hom(\Omega^{n-1},L) - \{0\})/\Gm$. 

\medskip

 Now we are able to prove the following result.

\begin{Pp} 
\label{Pp_majus'}
1) If Conjecture~\ref{Con_1} is true then 
$\DD(\Aut_{\bar E})\iso \Aut_{\bar E^*}$ canonically.

\noindent
2) If Conjecture~\ref{Con_1} is true for $\bar R=\cO$ then it is 
true in full generality. Moreover, there exists a perverse
sheaf $\Autb_{\bar E}\in\Perv_{fl}(\Bunb_n,\bar R)_{\sigma}$ together with
an isomorphism 
\begin{equation}
\label{iso_Autb}
\Aut_{\bar E}\iso\gr^*\Autb_{\bar E}[-1](-\frac{1}{2})
\end{equation}
By this condition $\Autb_{\bar E}$ is defined up to a canonical isomorphism. 
The formation of $\Autb_{\bar E}$ commutes with extension of scalars 
(\ref{ext_scalars}). We also have canonically 
\begin{equation}
\label{iso_multb}
\multb^*\Autb_{\bar E}\iso\Autb_{\bar E}\boxtimes\wedge^n \bar E
\end{equation} 
\end{Pp} 

\begin{Prf}
1) Let $\Bun_n^{est}\subset\Bun_n$ be the open substack consisting of $L$
such that $\deg L>n^2(g-1)$ and $\Hom(L,\cL^{est})=0$. Let 
$_n\cM^{est}=\cup_{d\ge 0}\;{_n\cM^{est}_d}$. By Conjecture~\ref{Con_1}, 
over $_n\cM^{est}$ there is a canonical isomorphism 
\begin{equation}
\label{iso_1}
\DD(_n\cK_{\bar E})\iso{_n\cK_{\bar E^*}}  
\end{equation}
Over $\Bun_n^{est}\cap\Bun_n^d$ the map 
$\varrho:{_n\cM_d}\to\Bun_n^d$ 
is a vector bundle of rank $d-n^2(g-1)$ with removed zero section. Since
the preimage of $\Bun_n^{est}\cap\Bun_n^d$ under this map is contained in
$_n\cM^{est}_d$, the isomorphism (\ref{iso_1}) descends to give an isomorphism
over $\Bun_n^{est}$
\begin{equation}
\label{iso_2}
\DD(\Aut_{\bar E})\iso\Aut_{\bar E^*}
\end{equation}
By (\cite{FGV2}, 1.5), Hecke property of $\Aut_{\bar E}$ yields an isomorphism
\begin{equation}
\label{iso_3}
\mult^*\Aut_{\bar E}\iso \Aut_{\bar E}\boxtimes\wedge^n\bar E
\end{equation}
For any open substack
of finite type $U\subset\Bun_n$ there exists an integer $d'$ such that for any
$x\in X$ the morphism $\mult_{d'x}:\Bun_n\to\Bun_n$ sending $L$ to $L(d'x)$
maps $U$ isomorphically onto a substack of $\Bun_n^{est}$. We get
\begin{multline*}
\DD(\Aut_{\bar E})\mid_U\iso \;\DD((\wedge^n\bar E_x)^{\otimes\,  -d'}\otimes
\mult^*_{d'x}\Aut_{\bar E})\;\iso\\
(\wedge^n\bar E^*_x)^{\otimes\, -d'}\otimes
\mult^*_{d'x}(\DD\Aut_{\bar E}) \;\iso\;
(\wedge^n\bar E^*_x)^{\otimes -d'}\otimes\mult^*_{d'x}\Aut_{\bar E^*}\iso
\Aut_{\bar E^*}\mid_U
\end{multline*}
According to (\ref{iso_3}), this gives a well-defined isomorphism (\ref{iso_2}) 
over the entire $\Bun_n$, which coincides with the old one over $\Bun_n^{est}$. 

\smallskip
\noindent
2) {\bf Step 1.} Applying the conservative functor $\EE\Lotimes_{\bar R_{\sigma}}$, 
one checks that $\zeta_!({_n\cF^d_{\bar E,\psi}})\to 
\zeta_*({_n\cF^d_{\bar E,\psi}})$ is an isomorphism over 
$_n\cM_d^{est}$. This yields the isomorphism (\ref{iso_1})
over $_n\cM^{est}_d$. 

 Since both $_n\cK^d_{\bar E}\Lotimes_{\bar R_{\sigma}}\EE$ and
$(\DD\, {_n\cK^d_{\bar E}})\Lotimes_{\bar R_{\sigma}}\EE$ 
restricted to $_n\cM^{est}_d$ are irreducible perverse sheaves, by 2) 
of Lemma~\ref{Lm_for_Pp_majus'}, the restriction of $_n\cK^d_{\bar E}$ to 
$_n\cM^{est}_d$ is an object of $\Perv_{fl}({_n\cM^{est}_d},\bar R)_{\sigma}$, 
which is the Goresky-MacPherson extension from any nonempty open substack of 
$_n\cM^{est}_d$. 

\smallskip
\Step 2 As was explained in (\cite{FGV2}, 7.6), we have a notion of a
Hecke eigensheaf on $\Bun_n^{est}$. Indeed, for the diagram
(\ref{diag_Hecke}) we have 
$$
(\gq\times\supp)^{-1}(\Bun_n^{est}\times X)\subset \gp^{-1}(\Bun_n^{est})
$$ 
Define the functor $\D^b_c(\Bun_n^{est},\bar R)_{\sigma}\to
\D^b_c(\Bun_n^{est}\times X,\bar R)_{\sigma}$ by formula
(\ref{hecke_functor}) and denote it again by $\H^1_n$. Considering its iterations
$(\H^1_n)^{\boxtimes d}$, one can repeat Definition~\ref{Def_Hecke_eigensheaf} in
this context.  

 Let $\Bunb_n^{est}$ be the image of $\Bun_n^{est}$ under $\gr:\Bun_n\to\Bunb_n$. 

\begin{Lm}
1) There exists a perverse sheaf $\Autb_{\bar E}^{est}\in\Perv_{fl}(\Bunb_n^{est},
\bar R)_{\sigma}$ and for each $d\ge 0$ an isomorphism
$$
\varrho^*\gr^*\Autb_{\bar E}^{est}[d-1-n^2(g-1)](\frac{d-1-n^2(g-1)}{2})\iso 
{_n\cK^d_{\bar E}}
$$
over $\varrho^{-1}(\Bun_n^d\cap\Bun_n^{est})$. These data are defined up to a 
canonical isomorphism. Besides, over each connected component of 
$\Bunb_n^{est}$, $\Autb_{\bar E}^{est}$ is the Goresky-MacPherson 
extension from any nonempty open substack, and the formation of 
$\Autb_{\bar E}^{est}$ commutes with extension 
of scalars (\ref{ext_scalars}). 

\medskip
\noindent
2) Set $\Aut_{\bar E}^{est}$ to be $\gr^*\Autb_{\bar E}^{est}[-1](-\frac{1}{2})$ 
over $\Bun_n^{est}$. Then $\Aut_{\bar E}^{est}$ has a unique structure of a 
Hecke eigensheaf with respect to $\bar E$, which is compatible with the Hecke 
property of $_n\cK_{\bar E}$. 
\end{Lm}
\begin{Prf}
Let $d>n^2(g-1)$ be such that $\Bun_n^d\cap\Bun_n^{est}$ is nonempty. 
Let $U\subset \Bun_n^d\cap\Bun_n^{est}$ be a nonempty open substack
such that $\Aut_{E_0}$ is an (appropriately shifted) smooth $\EE$-sheaf over $U$.
Let $U'=\varrho^{-1}(U)$. So, over $U'$, 
$_n\cK^d_{\bar E}$ is a smooth $\bar R_{\sigma}$-sheaf (appropriately shifted).  

 Recall that $U'\to U$ is a vector
bundle of rank $d-n^2(g-1)$ with removed zero section. Let $\bar U$ be the image
of $U$ under $\gr:\Bun_n\to\Bunb_n$. Then $U'\to\bar U$ is a projectivization
of a vector bundle. By 1) of Lemma~\ref{Lm_for_Pp_majus'}, 
there exists an (appropriately shifted)
smooth $\bar R_{\sigma}$-sheaf $V_0$ on $\bar U$ and an isomorphism in 
$\D^b_c(U',\bar R)_{\sigma}$
$$
\varrho^*\gr^*V_0[d-1-n^2(g-1)](\frac{d-1-n^2(g-1)}{2})\iso {_n\cK^d_{\bar E}}
$$ 
 
 Define the restriction of $\Autb_{\bar E}^{est}$ to $\Bunb_n^d\cap\Bunb_n^{est}$ 
as the Goresky-MacPherson extension of $V_0$ from $\bar U$. Since 
$\varrho^{-1}(\Bun_n^{est})\subset {_n\cM^{est}}$, our first assertion follows.
All the other assertions follow from the fact that the restriction
$\varrho^{-1}(\Bun_n^{est})\to\Bun_n^{est}$ of $\varrho$ is smooth and surjective 
with connected fibres. 
\end{Prf}
 
\medskip
\Step 3  As in (\cite{FGV2}, 1.5), Hecke property of $\Aut_{\bar E}^{est}$ yields
the isomorphism (\ref{iso_3}) over $\Bun_n^{est}\times X$. This isomorphism descends to
give the isomorphism (\ref{iso_multb}) over $\Bunb_n^{est}\times X$.
Then one extends $\Autb_{\bar E}^{est}$ to the entire $\Bunb_n$
as follows. 

 For any open substack
of finite type $U\subset\Bunb_n$ there exists an integer $d'$ such that for any
$x\in X$ the morphism $\multb_{d'x}:\Bunb_n\to\Bunb_n$ sending $L$ to $L(d'x)$
maps $U$ isomorphically onto a substack of $\Bunb_n^{est}$. Set $\Autb_{\bar E}\mid_{U}$
to be 
$$
(\wedge^n \bar E_x)^{\otimes -d'}
\otimes\multb^*_{d'x}\Autb_{\bar E}^{est}
$$
This gives a well-defined perverse sheaf $\Autb_{\bar E}\in\Perv_{fl}(\Bunb_n)$ together
with the isomorphism (\ref{iso_multb}) over the entire $\Bunb_n\times X$. One concludes
the proof as in (\select{loc.cit.}, 7.8 and 7.9). \\
\end{Prf}(Proposition~\ref{Pp_majus'})

\medskip
\begin{Rem}
\label{Rem_PGL}
1) Denote by $\Bun_{\PGL_n}$ the moduli stack of $\PGL_n$-bundles on $X$. There
exists a morphism $\tilde\gr: \Bunb_n\to\Bun_{\PGL_n}$ such that the
composition
$\Bun_n\toup{\gr}\Bunb_n\toup{\tilde\gr}\Bun_{\PGL_n}$ is the
canonical map $\Bun_n\to \Bun_{\PGL_n}$. The morphism $\tilde\gr$
is representable, smooth and separated. Let
$t_X:\Bun^d_n\times\Pic^0 X\to\Bun^d_n$ be the map that sends $(L,\cA)$ to
$L\otimes\cA$. We have a map $\u{t}_X:\Bunb^d_n\times\uPic^0
X\to\Bunb^d_n$ such that the diagram
$$
\begin{array}{ccc}
\Bun^d_n\times\Pic^0 X & \toup{t_X} & \Bun^d_n\\
\downarrow && \downarrow\\
\Bunb^d_n\times\uPic^0 X & \toup{\u{t}_X} & \Bunb^d_n
\end{array}
$$
commutes, and the following two squares are cartesian
$$
\begin{array}{ccccccccc}
\Bun^d_n & \to & \Bun_{\PGL_n} &&&
\Bunb^d_n & \toup{\tilde\gr} & \Bun_{\PGL_n}\\
\uparrow\lefteqn{\scriptstyle{t_X}} && \uparrow &&&
\uparrow\lefteqn{\scriptstyle{\u{t}_X}} &&
\uparrow\lefteqn{\scriptstyle{\tilde\gr}}\\
\Bun^d_n\times\Pic^0 X & \toup{\pr_1} & \Bun^d_n &&&
\Bunb^d_n\times\uPic^0 X & \toup{\pr_1} & \Bunb^d_n
\end{array}
$$
2) Let $\Lambda$ be a noetherian ring such that the characteristic of 
$\Lambda$ is invertible in $k$. Then $\Bunb_n$ is a Bernstein-Lunts 
stack with respect to $\Lambda$ in the sense of (\cite{LMb}, 18.7.4). 
 
 Indeed, if $\cX_1\to\cX_2$ is
a representable separated morphism of algebraic stacks, and $\cX_2$ is a
Bernstein-Lunts stack then $\cX_1$ is also. Apply this for $\tilde\gr$.
The stack $\Bun_{\PGL_n}$ is a Bernstein-Lunts stack, because it is of the form
$M/G$, where $M$ is a separated algebraic space with an action of an affine
algebraic group $G$ (\cite{LMb}, 18.7.5). 
\end{Rem}

\section{Deformations of local systems and cohomology of $\HOM(\bar E_1,\bar E_2)$}
\label{sect_univ_deform}

3.1 Let $\EE$ be a finite extension field of $\Ql$. Fix a smooth
$\EE$-sheaf $E_0$ on $X$ of rank $m$. First, we recall  the structure of the
universal deformation of $E_0$ over $\EE$.  This construction is standard
(cf. \cite{Schle} for the definition of pro-representability, etc.). 

 Let $\eta\in X$ be the generic point of $X$ and $\bar\eta\to\eta\to X$
be a geometric point over $\eta$. Set $G=\pi_1(X,\bar\eta)$. Let $A\in
\Ob(\cC_{\,\EE})$. Recall that the
functor that sends $E$ to $E_{\bar\eta}$ is an
equivalence between the category of smooth $A$-sheaves 
of rank $m$ on $X$ 
and the category of pairs $(E,\rho)$, where $E$ is a free $A$-module of rank
$m$ and $\rho: G\to \Aut_A E$ is a representation 
continuous in the $\ell$-adic topology. 

\begin{Def}
An \select{$A$-deformation of $E_0$} is a pair
$(E, h)$, where $E$ is a smooth $A$-sheaf on $X$ of rank $m$ 
and $h: E\otimes_A \EE\, \iso E_0$ is an isomorphism 
of $\EE$-sheaves on $X$.
\end{Def}

Define the functor $F_{E_0}:\cC_{\,\EE} \to \Sets$ by
$F_{E_0}(A)=$ the  set of isomorphism classes of $A$-deformations of $E_0$.

\begin{Pp} 
\label{pro-repr} 
If $\End(E_0)=\EE$ then  
$F_{E_0}$ is pro-representable by a pro-pair $(R,E)$, where $R$ is
(non canonically) isomorphic to the ring of formal power series over 
$\EE$ in $2+(2g-2)m^2$ variables.  
Let $\gm\subset R$ be the maximal ideal of $R$. The $\EE$-dual of
$\gm/\gm^2$ is canonically identified with 
$\H^1(X,\END E_0)$. If, in addition, $\EE\subset \EE'$ is a finite extension
field and $\End(E_0\otimes_{\,\EE}\,\EE')=\EE'$ then the pro-pair
$(R\otimes_{\,\EE}\,\EE' \, , E\otimes_{\,\EE}\EE')$ pro-represents
the functor $F_{E_0\otimes_{\,\EE}\;\EE'}$. 
\end{Pp}

\begin{Lm} 
\label{Schle-lemma}
Suppose that $\End(E_0)=\EE$\\
1) If $E$ is an $A$-deformation of $E_0$ then $\End(E)=A$.\\
2) Let $A'\to A, A''\to A$ be two morphisms in $\cC_{\,\EE}$. Suppose
that $A''\to A$ is surjective then the natural morphism
$
F_{E_0}(A'\times_A A'')\to F_{E_0}(A')\times_{F_{E_0}(A)} F_{E_0}(A'')
$
is a bijection.
\end{Lm}
\begin{Prf}
2) The surjectivity is easy. To prove the injectivity use point 1) and
Corollary 3.6, p.217 of \cite{Schle}.
\end{Prf}

\begin{Lm}
\label{Lm_for_formal_smoothness}
Let $A'\to A$ be a surjection in $\cC_{\,\EE}$ , whose
kernel is a 1-dimensional $\EE$-vector space. 
Let $V'$ be a free $A'$-module of rank
$m$. Put $V=V'\otimes_{A'} A$. Then the natural map
$$
\GL(V')\to
\PGL(V')\times_{\PGL(V)}\GL(V)\times_{\GL(\det V)}\GL(\det V')
$$
is an isomorphism of groups. 
\QED
\end{Lm}

\smallskip

\begin{Prf}\select{of Proposition \ref{pro-repr}}\\
Consider the ring $\EE[\varepsilon]/(\varepsilon^2)$ of dual numbers.
The groupo\"\i d of $\EE[\varepsilon]/(\varepsilon^2)$-deformations
of $E_0$ is naturally equivalent to the category of extensions
$0\to E_0\to ?\to E_0\to 0$ on $X$. So, the tangent space
to $F_{E_0}$ is identified with $\Ext^1_X(E_0,E_0)=\H^1(X,\END E_0)$.
Now combining Lemma \ref{Schle-lemma} with Theorem 2.11 of
\cite{Schle} we get the pro-representability of $F_{E_0}$ by a pro-pair
$(R,E)$. 

 Let us show that the morphism of functors associating to an $A$-deformation
$E$ of $E_0$ the $A$-deformation 
$\det E$ of $\det E_0$  is a formally smooth morphism from the
universal deformation of $E_0$ to the universal deformation of 
$\det E_0$. Let $A'\to A$, $V'$ and $V$ be as in
Lemma~\ref{Lm_for_formal_smoothness}. Suppose that $V$ is equipped with a
structure of an $A$-deformation of $E_0$. Let $\rho:G\to \Aut_A V$ be the
corresponding representation of $G$.  Since $\H^0(X,\END_0 E_0)=0$, we get
$\H^2(X,\END_0 E_0)=0$ (we write $\END_0 E$ for the sheaf of traceless 
endomorphisms of $E$). It follows that the 
corresponding representation of $G$ in $\PGL(V)$ can be lifted to a
representation $\rho':G\to \PGL(V')$. Now our assertion follows from
Lemma~\ref{Lm_for_formal_smoothness}.

 The universal deformation of $\det E_0$ is formally
smooth,  because it  is isomorphic to the universal deformation of the
trivial  1-dimensional local system, which is an infinitesimal formal
$\EE$-group (cf. SGA3, t.1,$\on{VII_B}$, 3.3).
So, $R$ is formally smooth, i.e., by (\cite{Schle}, 2.5), is isomorphic to the ring of formal power series over 
$\EE$. 

 Since $\chi(\END_0 E_0)=(2-2g)(m^2-1)$, we have 
$\dim \H^1(X,\END E_0)=(2g-2)m^2+2$.

 If $\EE\subset \EE'$ is a finite exension with
$\End(E_0\otimes_{\EE}\EE')=\EE'$ then $(R\otimes \EE' \, ,
E\otimes \EE')$ is a pro-pair for the functor 
$F_{E_0\otimes\EE'}:\cC_{\,\EE'}\to\Sets$, which defines a
morphism of functors $h_{R\otimes\EE'}\to F_{E_0\otimes\EE'}$,
where $h_{R\otimes\EE'}:\cC_{\,\EE'}\to\Sets$ denotes the
functor represented by $R\otimes\EE'$, that is, 
$h_{R\otimes\EE'}(B)=\Hom_{local\;\, \EE'-alg}(R\otimes\EE' \, , B)$. 
We must show that this is an isomorphism of functors. 
Since $F_{E_0\otimes\EE'}$ can
be represented by a ring of formal power series over $\EE'$, our assertion
follows from the fact that the  induced map on the tangent spaces is an
isomorphism. 
\end{Prf}

\medskip

 By definition, $E=(E_k)_{k\in\NN}$, $E_k\in F_{E_0}(R/\gm^k)$ are
such that the image of $E_{k+1}$ under $F_{E_0}(R/\gm^{k+1})\to 
F_{E_0}(R/\gm^k)$ is $E_k$. Fix a $R/\gm^k$-deformation of $E_0$
in the isomorphism class $E_k$ and denote it by the same 
symbol $E_k$. For each $k$ fix an isomorphism of 
$R/\gm^k$-deformations of $E_0$: $E_{k+1}
\otimes_{R/\gm^{k+1}} R/\gm^k \iso E_k$.
Then the projective system $(E_k)_{k\in \NN}$ is an object of $\Sh(X,R)$,
equipped with an isomorphism  $\alpha: E\otimes_R \EE \iso E_0$  of
$\EE$-sheaves on $X$. 
Notice that $R$ is defined up to a canonical isomorphism, whence the
$R$-sheaf $E$ is defined up to a non-canonical isomorphism.

\medskip
\noindent
3.2 Let $\cO\subset\EE$ be the ring of integers of $\EE$, $\kappa$ be the
residue field of $\cO$, and $\omega\in\cO$ be a uniformizing parameter. 
A smooth $\cO$-sheaf $E'_0$ on $X$ together with an
isomorphism $E'_0\otimes_{\cO}\EE\iso E_0$ can be viewed as a $G$-invariant
$\cO$-lattice in $(E_0)_{\bar\eta}$. Set $\bar E_0=E'_0\otimes_{\cO}\kappa$.
It is easy to see that, though $\bar E_0$ is not defined up to an
isomorphism by $E_0$, the image of $\bar E_0$ in the Grothendieck group
of the category of smooth $\kappa$-sheaves on $X$ is uniquely defined 
by $E_0$. 

 In this subsection we compare the universal deformation $(R,E)$ of $E_0$
and the universal deformation of $\bar E_0$ over $\cO$. 

 Let $\cC_{\cO}$ be the category of local Artin $\cO$-algebras with
residue field $\kappa$ (the morphisms are local homomorphisms of
$\cO$-algebras). For $A\in\Ob(\cC_{\cO})$ one defines a notion of an 
\select{$A$-deformation of $\bar E_0$} and a functor $F_{\bar E_0}:
\cC_{\cO}\to\Sets$ as in Sect. 3.1. The proof of the next result is
similar to that of Proposition~\ref{pro-repr}.

\begin{Pp}
\label{Pp_useful_however}
 If $\End(\bar E_0)=\kappa$ then $F_{\bar E_0}$ is
pro-representable by a pro-pair $(\bar R,\bar E)$, where $\bar R$ is 
a complete local noetherian $\cO$-algebra with residue field $\kappa$. If
$\bar m$ is the maximal ideal of $\bar R$ then the $\kappa$-dual of $\bar
m/(\bar m^2,\omega)$ is canonically identified with $\H^1(X, \END \bar E_0)$.
\QED
\end{Pp}

As in Sect. 3.1, we may and will view $\bar E$ as an object of $\Sh(X,\bar
R)$ equiped with an isomorphism $\bar\alpha:\bar E\otimes_{\bar R}\kappa\iso
\bar E_0$.  (Notice that $\bar R$ is defined by $\bar E_0$ up to a
canonical isomorphism, whence $\bar E$ is defined up to 
a non canonical isomorphism.)

 Since $E'_0$ is an $\cO$-deformation of $\bar E_0$, it defines a local
homomorphism of $\cO$-algebras $\sigma:\bar R\to\cO$ such that there exists
an isomorphism $\beta: \bar E\otimes_{\bar R}\cO\iso E'_0$ in $\Sh(X,\cO)$
compatible with $\bar\alpha$ (so, $\beta$ is defined up to multiplication
by an element of $1+\omega\cO$). 
Let $\gm_{\sigma}$ denote the kernel of the
induced map $\bar R\otimes_{\cO}\EE\to \EE$. We denote by $\bar R_{\sigma}$
the localization of $\bar R\otimes_{\cO}\EE$ in $\gm_{\sigma}$, and by
$\hat R_{\sigma}$ the $\gm_{\sigma}$-adic completion of $\bar
R\otimes_{\cO}\EE$.

\begin{Rem} The rings $\bar R\otimes_{\cO}\EE$ (and, hence, $\bar
R_{\sigma}$ and $\hat R_{\sigma}$) are noetherian. Indeed, $\bar R$ is
isomorphic to a quotient of the ring of formal power series
$\cO[[t_1,\ldots, t_N]]$ for some $N$, so that $\bar R\otimes_{\cO}\EE$ is
isomorphic to a quotient of 
$\cO[[t_1,\ldots, t_N]]\otimes_{\cO}\EE$, which is noetherian (cf. \cite{GL},
Appendix A.2).
\end{Rem}

 The isomorphism $\beta$ induces on $\bar E\otimes_{\bar R} \hat
R_{\sigma}$ a structure of a $\hat R_{\sigma}$-deformation of $E_0$.
This defines a local homomorphism of $\EE$-algebras $\tau: R\to \hat
R_{\sigma}$ that does not depend on the choice of $\beta$.

\begin{Pp}
\label{Pp_bar_R_is_smooth} 
Assume that $\End(E_0)=\EE$ and
$\End(\bar E_0)=\kappa$. Then $\cO\to\bar R$ is formally smooth,
that is, $\bar R$ is (non canonically) isomorphic to the ring of formal 
power series over $\cO$ in $2+(2g-2)m^2$ variables. Besides, 
the natural map $\tau: R\to \hat R_{\sigma}$ is an isomorphism of 
$\EE$-algebras.   
\end{Pp}

\begin{Rem} One easily checks that if $\bar E_0$ is an irreducible
$\kappa$-sheaf then $E_0$ is an irreducible $\EE$-sheaf on $X$.
However, \select{the converse is not true}. 
\end{Rem}

 We start with the following observation. Let $A\in\cC_{\,\EE}$. Denote by
$\gn\subset A$ the maximal ideal. If $A'\subset A$ is an 
$\cO$-subalgebra, which is a finite $\cO$-module with
$A'\otimes_{\cO}\EE\iso A$, then $A'=\cO\oplus \gn'$, where $\gn'=A'\cap
\gn$ is an $\cO$-lattice in the $\EE$-vector space $\gn$. Conversely,
if $\gn'$ is an $\cO$-submodule of finite type in $\gn$ with
$\gn'\otimes_{\cO}\EE\iso\gn$ then the $\cO$-subalgebra $A'\subset A$
generated by $\gn'$ satisfies $A'\otimes_{\cO}\EE\iso A$, and $A'$ is a
finite $\cO$-module.  

\begin{Lm}
\label{Lm_submodules}
1) Let $A'\subset A$ be an $\cO$-subalgebra, which is a finite
$\cO$-module with $A'\otimes_{\cO}\EE\iso A$. Let $M$ be a free $A$-module
of rank $m$, $M'\subset M$ be an $A'$-submodule of finite type with
$M'\otimes_{\cO}\EE\iso M$. Then there exists an $\cO$-subalgebra 
$A'\subset A''\subset A$ such that $A''$ is a finite $\cO$-module, and
$A''M'$ is a free $A''$-module. \\
2) Let, in addition, $M'_0$ denote the image of $M'\to
M/\gn M$ then the image of $A''M'\to M/\gn M$ also equals $M'_0$, and
the latter map induces an isomorphism of $\cO$-modules
$A''M'\otimes_{A''}\cO\iso M'_0$. 
\end{Lm}

\begin{Prf}
Notice that $M'_0=M'/(\gn M\cap M')$ is a free $\cO$-submodule of rank
$m$  (an $\cO$-lattice) in the $m$-dimensional $\EE$-vector space $M/\gn
M$.  Pick $e_1,\ldots,e_m\in M'$ that define an $\cO$-basis in $M'_0$. 
Then $e_1,\ldots,e_m$ is an $A$-basis in $M$. 
Set $L=A'e_1+\ldots+A'e_m\subset M'$. 

 Any $x\in M'$ is written uniquely as $x=a_1e_1+\ldots+a_me_m$ with $a_i\in
A$. Reducing modulo $\gn$ we learn that every $a_i$ lies in $\cO\oplus\gn$.

 Let $u_1,\ldots, u_r$ generate $M'$ over $A'$. Write $u_i=\sum_j
a_{ij} e_j$ with $a_{ij}\in\cO\oplus\gn$. Let $A''$ be the $A'$-subalgebra
of $A$ generated by all $a_{ij}$. Then $M'\subset A''L$, so that
$A''M'=A''L$, and $A''L$ is a free $A''$-module. The second statement
is clear. 
\end{Prf}

\begin{Lm}
\label{Lm_entire_structure_on_smooth_sheaf}
1) Let $V$ be a smooth $A$-sheaf on $X$.
Then there exists an $\cO$-subalgebra $A'\subset A$, which is a finite
$\cO$-module with $A'\otimes_{\cO}\EE\iso A$, a smooth $A'$-sheaf
$V'$ on $X$, and an isomorphism $\nu: V'\otimes_{A'} A\iso V$ of
$A$-sheaves on $X$. \\
2) Given, in addition, a smooth $\cO$-sheaf $V'_0$
and an isomorphism $V'_0\otimes_{\cO}\EE\iso V\otimes_A \EE$,
one may chose $A', V',\nu$ such that the isomorphism
$(V'\otimes_{A'}\cO)\otimes_{\cO}\EE\to V'_0\otimes_{\cO}\EE$ induced by
$\nu$ is obtained by extension of scalars from an isomorphism
of $\cO$-sheaves $V'\otimes_{A'}\cO\iso V'_0$. 
\end{Lm}
\begin{Prf} 1) View $V$ as a free $A$-module with a continuous
representation of $G$. Pick a $G$-invariant $\cO$-lattice $V_1$ in
$V$. Pick any $\cO$-subalgebra $A'\subset A$, which is a
finite $\cO$-module with $A'\otimes_{\cO}\EE\iso A$, 
and apply Lemma~\ref{Lm_submodules} for
the $A'$-submodule $A'V_1$ of $V$.

\smallskip
\noindent
2) Recall that $\gn\subset A$ denotes the maximal ideal. Chose
$V_1$ with an additional property: the image of $V_1$ in $V/\gn V$
is the $\cO$-lattice $V'_0$. Then for any $\cO$-subalgebra 
$A''\subset A$, which is a finite $\cO$-module
with $A''\otimes_{\cO}\EE\iso A$, we still have the same property
for $V_1$ replaced by $A''V_1$. So, our assertion follows from ii) of
Lemma~\ref{Lm_submodules}. 
\end{Prf}

\medskip

\begin{Prf}\select{of Proposition~\ref{Pp_bar_R_is_smooth}}

\noindent
 Let $D'\to D$ be a surjection in $\cC_{\cO}$, whose kernel is a
principal ideal $(t)\subset D'$ with $t\,\gm_{D'}=0$, where $\gm_{D'}$
denotes the maximal ideal of $D'$. We must show that any local homomorphism
of $\cO$-algebras $\bar R\to D$ can be lifted to $\bar R\to D'$. 

 Pick $k>0$ with $\gm^k_D=0$, where $\gm_D\subset D$ is the maximal
ideal. Then $\gm^{k+1}_{D'}=0$. Set $A=R/\gm^{k+1}$ and denote by
$\gn\subset A$ the maximal ideal. By
Lemma~\ref{Lm_entire_structure_on_smooth_sheaf},  there is an
$\cO$-subalgebra $A'\subset A$, which is a finite $\cO$-module  with
$A'\otimes_{\cO}\EE\iso A$, a smooth $A'$-sheaf $V$ on $X$ and an 
isomorphism $\delta: V\otimes_{A'}A\iso E\otimes_R A$
of $A$-sheaves on $X$ such that the induced isomorphism
$$
(V\otimes_{A'}\cO)\otimes_{\cO}\EE\iso
E'_0\otimes_{\cO}\EE
$$ 
of $\EE$-sheaves is obtained by extension of scalars from an
isomorphism $\gamma: V\otimes_{A'}\cO\iso E'_0$ of $\cO$-sheaves on $X$.  
Further, $\gamma$ endows $V$ with a structure of an $A'$-deformation of 
$\bar E_0$. This defines a local homomorphism of $\cO$-algebras 
$\bar R\to A'$ such that $\bar E\otimes_{\bar R} A'$ and $V$ are 
isomorphic as $A'$-deformations of $\bar E_0$. 

 Notice that the composition $\bar R\to A'\to\cO$ coincides with $\sigma$.
Set $B=(\bar R\otimes_{\cO}\EE)/\gm_{\sigma}^{k+1}$, so that $\bar
R\to A'$ yields a morphism $B\to A$ in $\cC_{\,\EE}$. 

\begin{Lm} 
\label{Lm_two_homs_coincide}
The composition $R\toup{\tau} \hat R_{\sigma}\to B\to A$ 
coincides with the natural map $R\to R/\gm^{k+1}$. 
\end{Lm}
\begin{Prf} 
Let $f_i:R\to C$ be two local homomorphisms of $\EE$-algebras with
$C\in\Ob(\cC_{\,\EE})$.  If $h: E\otimes_{R,f_1} C\iso E\otimes_{R, f_2}
C$ is an isomorphism of $C$-sheaves on $X$ then, for a suitable $a\in\EE^*$,
$ah$ is an isomorphism of $C$-deformations of $E_0$, which implies 
$f_1=f_2$. In our case the corresponding $A$-sheaves on $X$ are 
isomorphic by definition. 
Indeed, by definition of $\tau:R\to \hat R_{\sigma}$, we have 
$(E\otimes_R \hat R_{\sigma})\otimes_{\hat R_{\sigma}} B
\iso \bar E\otimes_{\bar R} B$. Finally, 
$$
(\bar E\otimes_{\bar R}
B)\otimes_B A\iso (\bar E\otimes_{\bar R} A')\otimes_{A'} A\iso
V\otimes_{A'} A\toup{\delta} E\otimes_R R/\gm^{k+1}
$$
\end{Prf}

 By Lemma~\ref{Lm_two_homs_coincide}, the image $A_0$ of $\bar R\to A'$
generates $A$ as a $\EE$-vector space.
Replacing $A'$ by $A_0$ and $V$ by $\bar E\otimes_{\bar R} A_0$,
we may assume that $\bar R\to A'$ is surjective.

 Consider the composition $\bar R\to A'\to A\to A/\gn^2$. Its image 
is an $\cO$-subalgebra $\cO\oplus\tilde \gn$ of $A/\gn^2$ such that
$\tilde \gn$ is an $\cO$-lattice in $\gn/\gn^2$. Let $r$ denote
the dimension of $\gn/\gn^2$ over $\EE$. Pick $e_1,\ldots,e_r\in \bar R$
whose images in $A/\gn^2$ define an $\cO$-basis of $\tilde\gn$. 
Let $\rho: \cO[[t_1,\ldots,t_r]]\to \bar R$ be the continuous homomorphism
of $\cO$-algebras that takes $t_i$ to $e_i$. 

 We claim that $\rho$ is surjective. To see this it suffices to show that
the reduction $\kappa[[t_1,\ldots,t_r]]\to \bar R/\omega\bar R$ of $\rho$
is surjective. Now, $\bar R\to \cO\oplus\tilde\gn$ induces a $\kappa$-linear
surjective map $\bar\gm/(\bar\gm^2,\omega)\to \tilde\gn/\omega\tilde\gn$. 
The latter is an isomorphism, because the dimension of 
$\bar\gm/(\bar\gm^2,\omega)$ equals $r$ by
Proposition~ \ref{Pp_useful_however}. So, $\rho$ is surjective and induces
an isomorphism $\cO[[t_1,\ldots,t_r]]/(t_1,\ldots,t_r)^{k+1}\iso A'$.

 Clearly, $\bar R\to D$ factors as $\bar R\to A'\to D$ and
$A'\to D$ can be lifted to $A'\to D'$, so that $\cO\to\bar R$ is formally
smooth. By Lemma~\ref{Lm_two_homs_coincide}, $\tau: R\to \hat R_{\sigma}$
is injective. Since $\tau$ induces an isomorphism on tangent spaces, it is
an isomorphism. \\
\end{Prf} (Propoistion~ \ref{Pp_bar_R_is_smooth})

\medskip

\begin{Rem}
If we do not assume that $\End(\bar E_0)=\kappa$ then $F_{\bar E_0}$
has a hull $(\bar R,\bar E)$ in the sense of Schlessinger \cite{Schle}, 
where $\bar R$ is a complete local noetherian $\cO$-algebra (defined up to
a non canonical isomorphism). One still can define the maps
$\sigma:\bar R\to\cO$ and $\tau:R\to \hat R_{\sigma}$ as above (they are
no more unique) and show that $\tau$ is injective. 
We conjecture that this map $\tau: R\to\hat R_{\sigma}$ is formally smooth, that is,
$\hat R_{\sigma}$ is a ring of formal power series over $R$.  
\end{Rem}

\bigskip
\noindent
3.3 \select{Cohomology of $\HOM(\bar E_1,\bar E_2)$}

\smallskip\noindent
In the rest of Sect.~\ref{sect_univ_deform} we assume the 
conditions of Proposition~\ref{Pp_bar_R_is_smooth} satisfied.  

 Recall that
$(\bar R,\bar E)$ denotes the universal deformation of $\bar E_0$ over
$\cO$. Notice that $\bar R\Cotimes_{\cO}\bar R$ is isomorphic to the ring
of formal power series over $\cO$ in $2r$-variables, where $r=2+(2g-2)m^2$.
Put $\bar E_i=\bar E\otimes_{\bar R}(\bar R\Cotimes_{\cO}\bar R)$, where
the homomorphisms $\bar R\to \bar R\Cotimes\bar R$ correspond to two
projections $\Specf(\bar R\Cotimes\bar R)\to \Specf R$. So, $\HOM(\bar
E_1,\bar E_2)$ is a smooth $\bar R\Cotimes\bar R$-sheaf on $X$.
Consider the map 
$$
\HOM(\bar E_1,\bar E_2)\to
\HOM(\bar E_1,\bar E_2)\otimes_{\bar R\Cotimes \bar R} \bar R\iso
\END(\bar E)\toup{\tr} \bar R
$$ 
Applying the functor $\H^2(X,\, .\,)$, we
get a canonical map
$\H^2(X,\HOM(\bar E_1,\bar E_2))\to \H^2(X,\bar R)\iso \bar R(-1)$. 

\begin{Pp} 
\label{cohom.of_HOM}
1) $\RG(X,\HOM(\bar E_1,\bar E_2))$ is an object of 
$\D^b_{coh}(\bar R\Cotimes \bar R)$ that can be represented by a complex
$(V^0\to V^1\to V^2)$ of free $\bar R\Cotimes \bar R$-modules with
$\rk V^0=\rk V^2=1$, $\rk V^1=(2g-2)m^2+2$ such that the differential in
$V$ is zero modulo the maximal ideal of $\bar R\Cotimes \bar R$. 
The complex $V$ is defined up to a non canonical isomorphism 
of complexes. \\
2) The canonial map $\H^2(X,\HOM(\bar E_1, \bar E_2))\to \bar R(-1)$ 
is an isomorphism of $\bar R\Cotimes \bar R$-modules.
\end{Pp}

\begin{Prf}
1) By (\cite{E}, 6.3), $\RG(X,\HOM(\bar E_1,\bar E_2))$ lies in 
$\D^b_{coh}(\bar R\Cotimes \bar R)$, and we have 
$$
\RG(X, \END(\bar E_0))\,\iso\, \RG(X,\HOM(\bar E_1,\bar E_2))
\Lotimes_{\bar R\Cotimes\bar R} \kappa
$$
By Lemma~\ref{Lm_perfect_derived}, $\RG(X,\HOM(\bar E_1,\bar E_2))$ can
be represented by a perfect complex $V$ of $\bar R\Cotimes\bar R$-modules
whose differential is zero modulo the maximal ideal of $\bar R\Cotimes\bar
R$, and $V$ is defined up to a non canonical isomorphism. However, the
complex $V\otimes_{\bar R\Cotimes\bar R} \kappa$ is defined up to a
canonical isomorphism. More presicely, $V^i\otimes_{\bar R\Cotimes \bar R}
\kappa \iso \H^i(X,\END(\bar E_0))$ canonically for every $i$. Our first
assertion follows.

\smallskip
\noindent
2) By the projection formulae (\cite{GL}, ii) of Proposition A.1.5),
$$
\RG(X,\HOM(\bar E_1,\bar E_2)) \Lotimes_{\bar R\Cotimes\bar R} \bar R
\,\iso\, \RG(X, \END \bar E)\,\iso\,
\RG(X,\bar R)\oplus \RG(X, \END_0 \bar E)
$$ 
Since $\H^i(X, \END_0 \bar E)=0$ for $i\ne 1$, 
the differential in $V\otimes_{\bar R\Cotimes\bar R} \bar R$ vanishes.
Let $A\in \Ob \cC_{\cO}$. Let $\bar R\to A$ be a surjective local 
homomorphism of $\cO$-algebras, $I\subset A\otimes_{\cO}A$ be the ideal of
the diagonal, and $J\subset I$ be another ideal. 

The next assertion is an immedialte consequence of the universal property
of $(\bar R,\bar E)$.

\begin{Lm}
\label{universal}
If the images of
$\bar E_1\otimes_{\bar R\Cotimes \bar R} (A\otimes A)$ and
$\bar E_2\otimes_{\bar R\Cotimes \bar R} (A\otimes A)$ in
$F_{\bar E_0}(A\otimes A/J)$ coincide then $J=I$. \QED
\end{Lm}

\begin{Lm}
\label{vanishes}
Let $B\in \Ob\cC_{\cO}$ and $M_1,M_2$ be two non-isomorphic $B$-deformations
of $\bar E_0$. Then for any $f\in\H^0(X,\HOM(M_1,M_2))$
the map $f\otimes\id: M_1\otimes_B \kappa \to M_2\otimes_B \kappa$ vanishes.
\end{Lm}
\begin{Prf}
We have $f\otimes\id \in \End(\bar E_0)=\kappa$. If $f\otimes\id \ne 0$ then
$f\otimes\id$ is an isomorphism of $\kappa$-vector spaces. Therefore,
$f$ is an isomorphism of smooth $B$-sheaves on $X$. After a multiplication
by a  suitable element of $\kappa^*$, $f$ becomes an isomorphism of 
$B$-deformations of $\bar E_0$.
\end{Prf}

\medskip
Set $B=A\otimes A/J$. 
\begin{Lm}
\label{on_d_0} If the differential 
$d^0: V^0\otimes_{\bar R\Cotimes \bar R}
B \to V^1\otimes_{\bar R\Cotimes \bar R}B$ vanishes then $J=I$.
\end{Lm}
\begin{Prf}
Consider the $B$-deformations 
$M_i=E_i\otimes_{\bar R\Cotimes \bar R} B$ of $\bar E_0$  $(i=1,2)$. By our
assumption,
$\H^0(X,\HOM(M_1,M_2))\iso V^0\otimes_{\bar R\Cotimes \bar R} B$ is a
free $B$-module of rank 1. Suppose that $J\ne I$ then
$M_1$ and $M_2$ are non-isomorphic by Lemma~\ref{universal}.
Denote by $\gn$ the maximal ideal of $B$ and set $\Ann\gn=
\{b\in B\mid b\gn=0\}$. By Lemma~\ref{vanishes}, $\Ann\gn$ annihilates
$\H^0(X,\HOM(M_1,M_2))$. Since $\Ann\gn\ne 0$, we get a contradiction.
\end{Prf}

\medskip

 Consider the complex $V\otimes_{\bar R\Cotimes \bar R}(A\otimes A)$.
Combining Lemma~\ref{on_d_0} with the Poincar\'e duality, one proves
that the image of the differential $d^1: V^1\otimes_{\bar R\Cotimes \bar R}
(A\otimes A) \to V^2\otimes_{\bar R\Cotimes \bar R}(A\otimes A)$ is 
$I(V^2\otimes_{\bar R\Cotimes \bar R}(A\otimes A))$. 
In other words, the natural map
$$
\H^2(X, \HOM(\bar E_1,\bar E_2)\otimes_{\bar R\Cotimes \bar R}(A\otimes
A))\to A(-1)
$$ 
is an isomorphism. Passing to the limit we get the desired assertion.  \\
\end{Prf}(Proposition~\ref{cohom.of_HOM})

\bigskip
\noindent
3.4 \select{ Cohomology of $^x(\HOM(\bar E_1, \bar E_2))^{(d)}$}

\smallskip\noindent
Recall that we write $\Pic^d X$ for the Picard stack classifying 
invertible sheaves of degree $d$ on $X$. Let $\uPic^d X$ be the
corresponding Picard scheme of $X$, so that the natural map
$\gr:\Pic^d X\to \uPic^d X$ is a $\Gm$-gerbe. Chose a closed point $x\in
X$. It defines a section $\alpha_x:\uPic^d X\to \Pic^d X$ of $\gr$.
Namely, if one considers $\uPic^d X$ as the moduli scheme of pairs $(\cA,t)$,
where $\cA\in\Pic^d X$ and $t:\cA_x\iso k$ is a trivialization 
of the geometric fibre at $x$ then $\alpha_x$ sends $(\cA,t)$ to $\cA$. 

 Define the $\Gm$-torsor $\alpha'_x: {^xX^{(d)}}\to X^{(d)}$ by the
cartesian square
$$
\begin{array}{ccc}
^xX^{(d)} & \toup{\alpha'_x} & X^{(d)}\\
\downarrow && \downarrow\\
\uPic^d X & \toup{\alpha_x} & \Pic^d X,
\end{array}
$$
where the right vertical arrow sends a divisor $D\in X^{(d)}$ to $\cO_X(D)$.

 Recall that $\HOM(\bar E_1, \bar E_2)$ is a smooth $\bar
R\Cotimes\bar R$-sheaf on $X$, so that $(\HOM(\bar E_1, \bar
E_2))^{(d)}$ is a constructable $\bar R\Cotimes\bar R$-sheaf on $X^{(d)}$
(cf. Sect.~1.4). Let
$^x(\HOM(\bar E_1, \bar E_2))^{(d)}$ denote its inverse image to
$^xX^{(d)}$. In this subsection we prove the following result.
\begin{Pp}
\label{Pp_coh_of_xHOM}
For $d>0$ there is a canonical isomorphism of $\bar R\Cotimes \bar R$-modules
$$
\H^{2d+2-i}_c(^xX^{(d)}, {^x(\HOM(\bar E_1, \bar E_2))^{(d)}})\iso 
\begin{cases}
\bar R(-d-1), & \mbox{if } \;\, i=0\\
0, & \mbox{if } \;\, 0<i<d,
\end{cases}
$$
where the $\bar R\Cotimes \bar R$-module structure on $\bar R$ is given via
the diagonal mapping $\bar R\Cotimes \bar R\to \bar R$. 
\end{Pp}
 
 This will be done using Proposition~ \ref{cohom.of_HOM} and
Appendices~\ref{Some_linear_algebra} and 
\ref{Chern_classes}.

\medskip
\noindent
3.4.1 \  To prove Proposition~\ref{Pp_coh_of_xHOM} we need the following 
linear algebra lemma.

 Let $A$ be a (commutative) ring of characteristic 0. Consider a complex of
$A$-modules $M=(A\to M^{-1}\toup{\eth} A)$, where $M^{-1}$ is a free $A$-module of
rank $r$. (So, $M^{-2}=M^0=A$  and $M^i=0$ for $i>0$ and  $i<-2$).
Suppose that there exists a basis $e_1,\ldots,e_r\in M^{-1}$ such that 
$\eth(e_1),\ldots,\eth(e_r)$ is a regular sequence for $A$. Denote by $I$ the image of
$\eth$. Let $\xi_1: M\to M[2]$ be a morphism of complexes such that the induced
map $M^{-2}\to M^0$ is an isomorphism. Define the morphism 
$\xi: \otimesl_{i=1}^d M\to (\otimesl_{i=1}^d M)[2]$
as $\xi_1\otimes\id\otimes\dots\otimes\id+\ldots+
\id\otimes\dots\otimes\id\otimes\xi_1$.
Then there is a (unique) morphism $\xi_d:\Sym^d(M)  \to
\Sym^d(M)[2]$ such that the diagram commutes
$$
\begin{array}{ccc}
\otimesl_{i=1}^d M & \toup{\xi} & \otimesl_{i=1}^d M[2]\\
\cup && \cup\\
\Sym^d(M) & \toup{\xi_d} & \Sym^d(M)[2]
\end{array}
$$
Notice that $\Sym^d(M)$ is a bounded complex of free $A$-modules of
finite type. 
\begin{lm}
\label{Lm_A.1}
Define the object $K\in \D_{\parf}(A)$ from the distinguished triangle
$K\to \Sym^d(M)\toup{\xi_d}\Sym^d(M)[2]$. Then $\H^0(K)\iso A/I$ and $\H^i(K)=0$
for $-d<i<0$ and for $i>0$. 
\end{lm}
 The proof is given in Appendix~\ref{Some_linear_algebra}.

\medskip\noindent
3.4.2 \ Denote by $Y_x: X^{(d-1)}\hook{} X^{(d)}$ the closed immersion 
that sends $D$ to
$D+x$. We consider $Y_x$ as a divisor on $X^{(d)}$ and write sometimes 
$Y_x\hook{} X^{(d)}$ for the same closed subscheme. Denote by $^{\p}Y_x$ 
the inverse image of $Y_x$ under $\sym:X^d\to X^{(d)}$. (In other words, 
the closed immersion $^{\p}Y_x\hook{} X^d$ is obtained from 
$Y_x\hook{} X^{(d)}$ by the base change $\sym:X^d\to X^{(d)}$). Denote by 
$^{\p}Y^i_x$ the inverse image of $x$ under $\pr_i: X^d\to X$. So, 
$^{\p}Y^i_x$ and $^{\p}Y_x$ are divisors on $X^d$, and we have 
$^{\p}Y_x={^{\p}Y^1_x}+\ldots+{^{\p}Y^d_x}$. 

 Consider the invertible sheaf $\cO(Y_x)$ on $X^{(d)}$. 
\begin{Lm}
$^xX^{(d)}$ is naturally isomorphic to the total space of $\cO(Y_x)$ with
removed zero section. 
\end{Lm}
\begin{Prf}
Denote by $Y^{univ}\hook{} X^{(d)}\times X$ the universal divisor. 
Clearly, the inverse image of $Y^{univ}$ under the closed immersion 
$X^{(d)}\times x\hook{} X^{(d)}\times X$ is the divisor $Y_x\times x$ on 
$X^{(d)}\times x$ with some multiplicity $r>0$. It is enough to show that 
$r=1$, i.e., to show that the following square is cartesian
$$
\begin{array}{ccc}
Y^{univ} & \hook{} & X^{(d)}\times X\\
\uparrow && \uparrow\\
X^{(d-1)}\times x & \hook{Y_x\times\id} & X^{(d)}\times x
\end{array}
$$
To do so, denote by $^{\p}Y^{univ}$ the inverse image of $Y^{univ}$ under 
$X^d\times X\toup{\sym\times \id} X^{(d)}\times X$. Since the inverse 
image of $^{\p}Y^{univ}$ under the closed immersion $X^d\times x\hook{}
X^d\times X$ is $^{\p}Y_x\times x$ with multiplicity one, our assertion 
follows. 
\end{Prf}

\medskip
\begin{Prf}\select{of Proposition~\ref{Pp_coh_of_xHOM}}\\
Let for brevity $\cR=\bar R\Cotimes_{\cO}\bar R$. 
By Lemma~\ref{prop_C1}, we have a distinguished triangle 
$$
(\alpha'_x)_!\cR\to \cR(-1)[-2]\toup{c}\cR
$$ 
in $\D^b_c(X^{(d)},\cR)$, 
where $c\in \H^2(X^{(d)},\cR(1))$ is the Chern class of
$\cO(Y_x)$. By K\"unneth's formulae,
$$
\H^2(X^d,\cR)=\mathop{\oplus}\limits_{i_1+\dots+i_d=2}
\H^{i_1}(X,\cR)\otimes\dots\otimes
\H^{i_d}(X,\cR)
$$
and $\H^2(X^{(d)},\cR)=\H^2(X^d,\cR)^{S_d}$. Denote by $c'$ the image of
$c$ in $\H^2(X^d,\cR(1))$. The construction of the Chern class 
is functorial, so that $c'$ is the Chern class of $\cO(^{\p}Y_x)$. 
Since $^{\p}Y_x={^{\p}Y^1_x}+\ldots+{^{\p}Y^d_x}$, we get 
\begin{multline*}
c'=c_1\otimes
1\otimes\dots\otimes 1+\ldots+1\otimes\dots\otimes 1\otimes c_1 \in \\
\H^2(X,\cR(1))\otimes\H^0(X,\cR)\otimes\dots\otimes\H^0(X,\cR)\oplus\ldots
\\
\oplus\H^0(X,\cR)\otimes\dots\otimes\H^0(X,\cR)\otimes\H^2(X,\cR(1))
\subset \H^2(X^d,\cR(1)),
\end{multline*}
where $c_1\in\H^2(X,\cR(1))$ is the Chern class of the invertible sheaf
$\cO(x)$ on $X$. Since $\cO(x)$ is of degree 1, we have $c_1\ne 0$. 
On $X^{(d)}$ we get a distinguished triangle 
$$
(\alpha'_x)_!{^x(\HOM(\bar E_1, \bar E_2))^{(d)}}(1)[2]\to
(\HOM(\bar E_1, \bar E_2))^{(d)}\toup{c} 
(\HOM(\bar E_1, \bar E_2))^{(d)}(1)[2]
$$
Let $\xi_d$ be the morphism obtained from
$(\HOM(\bar E_1,\bar E_2))^{(d)}\toup{c} (\HOM(\bar
E_1,\bar E_2))^{(d)}(1)[2]$  by applying the functor $\RG(X^{(d)}, \;\,
.\;\,)$. We get the distinguished triangle in
$\D^b_{coh}(\cR)$
\begin{multline*}
\RG_c(^xX^{(d)},{^x(\HOM(\bar E_1,\bar E_2))^{(d)}})(1)[2]\to
\RG(X^{(d)},(\HOM(\bar E_1,\bar E_2))^{(d)})
\toup{\xi_d} \\
\RG(X^{(d)},(\HOM(\bar E_1, \bar E_2))^{(d)})(1)[2]
\end{multline*}
Since $\sym_!(\HOM(\bar E_1,\bar E_2))^{\boxtimes\, d}$ is a direct sum 
over the irreducible representations of $S_d$, the same holds for 
$\RG(X^d, (\HOM(\bar E_1,\bar E_2))^{\boxtimes\,
d})=\otimesl_{i=1}^d\RG(X,\HOM(\bar E_1,\bar E_2))$, and we have  
naturally
$$
\RG(X^{(d)},(\HOM(\bar E_1,\bar E_2))^{(d)})\iso 
(\otimesl_{i=1}^d\RG(X,\HOM(\bar E_1,\bar E_2))\;)^{S_d}
$$
Denote also by $\xi: \otimesl_{i=1}^d\RG(X,\HOM(\bar E_1,\bar E_2)) \to
\otimesl_{i=1}^d\RG(X,\HOM(\bar E_1,\bar E_2))(1)[2]$ the morphism obtained
from $\sym_!(\HOM(\bar E_1,\bar E_2))^{\boxtimes\,
d}\toup{c}\sym_!(\HOM(\bar E_1,\bar E_2))^{\boxtimes\, d}(1)[2]$ by
applying the functor $\RG(X^{(d)}, \;\, .\;\,)$.

 The map $\xi$ is a cup product 
by an element $c\in \H^2(X^{(d)},\cR(1))$. Replacing
the cup product on $X^{(d)}$ by that on $X^d$, we get
$$
\xi=\xi_1\otimes\id\otimes\dots\otimes\id+\ldots+
\id\otimes\dots\otimes\id\otimes\xi_1
$$ 

 Pick a perfect complex $M$ of $\cR$-modules that represents
$\RG(X,\HOM(E_1,E_2))$. We assume that $M$ is chosen as in 1) of 
Proposition~\ref{cohom.of_HOM}. Pick a morphism 
$\tilde\xi_1:M\to M(1)[2]$ that
represents $\xi_1$ in $\D_{parf}(\cR)$ (so, $\tilde\xi_1$ is
defined up to a homotopy). Since $c_1\ne 0$, it follows that 
$\tilde\xi_1$ is given by the diagram
$$
\begin{array}{ccc}
& M^0 & \!\!\!\!\to M^1\to M^2\\
& \downarrow\\
M^0(1)\to M^1(1)\to\!\! & M^2(1),
\end{array}
$$
where the vertical arrow is an isomorphism of $\cR$-modules.

 Let $r=2+(2g-2)m^2$. By Proposition~\ref{cohom.of_HOM}, $M^1$ is a free
$\cR$-module of rank $r$. Further, any $r$ elements
in $\cR$, which generate the ideal of the diagonal, form a regular
sequence in $\cR$. So, combining 2) of Proposition~\ref{cohom.of_HOM} 
with Lemma~\ref{Lm_A.1} we get the desired assertion.
\end{Prf}

\section{Main Global Theorem}

4.1 Let $\EE$ be a finite extension of $\Ql$ that contains the group of 
$p$-th roots of unity, $\cO\subset\EE$ be its ring of integers, and 
$\kappa$ be the residue field of $\cO$. Let $E_0$ be a smooth $\EE$-sheaf 
of rank $n$ on $X$ such that $E_0\otimes_{\EE}\Qlb$ is irreducible (this, in
particular, implies $\End(E_0)=\EE$). 

 Choose a smooth $\cO$-sheaf $E'_0$ on $X$
together with an isomorphism $E'_0\otimes_{\cO}\EE\iso E_0$ and let
$\bar E_0=E'_0\otimes_{\cO}\kappa$. The local system $\bar E_0$ may not be
irreducible. We impose an additional assumption
$\End(\bar E_0\otimes_{\kappa}\bar\kappa)=\bar\kappa$, where $\bar\kappa$
is an algebraic closure of $\kappa$. This automatically implies
$\End(\bar E_0)=\kappa$.

  Let $(\bar R,\bar E)$ be the universal deformation of $\bar E_0$ over $\cO$.
It is equiped with a local homomorphism of $\cO$-algebras
$\sigma:\bar R\to\cO$ (cf. Sect. 3.2). 

 Put $\bar E_i=\bar E\otimes_{\bar R}
(\bar R\Cotimes_{\cO}\bar R)\;$ $(i=1,2)$, where the homomorphisms $\bar R\to
\bar R\Cotimes\bar R$ correspond to two projections 
$\Specf(\bar R\Cotimes\bar R)\to\Specf(\bar R)$. These are smooth 
$\bar R\Cotimes\bar R$-sheaves on $X$ of rank $n$. 

 By abuse of notation, the composition of the diagonal map 
$\bar R\Cotimes\bar R\to \bar R$ with $\sigma:\bar R\to\cO$ will also 
be denoted by $\sigma$. So, we have the categories $\D^b_c(\; .\; ,\bar
R\Cotimes\bar R)_{\sigma}\;$ (cf. Sect. 1.4). 
 
 Our main result is the following.
 
\begin{MGTh} Suppose that Conjecture~\ref{Con_1} is true. Then for any integer
$d$ there is a canonical isomorphism in $\D^b_c(\Spec k, \,
\bar R\Cotimes\bar R)_{\sigma}$
$$
\RG_c(\Bunb^d_n, \; \Autb_{\bar E_1^*}^d
\Lotimes_{(\bar R\Cotimes\bar R)_{\sigma}}\Autb_{\bar E_2}^d)\,\iso\,
\bar R
$$
where $\bar R$ is considered as a $\bar R\Cotimes\bar R$-module via 
the diagonal map $\bar R\Cotimes\bar R\to \bar R$.
\end{MGTh} 

\begin{Rems}
i) As in Sect. 1.4, let $(\bar R\Cotimes\bar R)_{\sigma}$ 
denote the localization of $\bar R\Cotimes\bar R$ in the 
multiplicative system $\{x\in \bar R\Cotimes\bar R\mid \sigma(x)\ne 0\}$.
Notice that $\bar R_{\sigma}$ is the localization of $\bar R$ in
$\{x\in\bar R\mid \sigma(x)\ne 0\}$. By Proposition~\ref{Pp_majus'}, 
$$
\Autb_{\bar E_i}\iso\Autb_{\bar E}\Lotimes_{\bar R_{\sigma}}(\bar R\Cotimes\bar
R)_{\sigma}
$$ 
is an object of $\Perv_{fl}(\Bunb_n,\bar R\Cotimes\bar R)_{\sigma}$. By
Remark~\ref{Rem_strange}, the core of the natural $t$-structure on 
$\D^b_c(\Spec k, \, \bar R\Cotimes\bar R)_{\sigma}$ is a full subcategory
of the category of $(\bar R\Cotimes\bar R)_{\sigma}$-modules. Since
$$
\bar R\otimes_{\bar R\Cotimes\bar R} 
(\bar R\Cotimes\bar R)_{\sigma}\iso \bar R_{\sigma},
$$ 
Main Global Theorem can be reformulated as follows: 
\select{for any integers $i,d$ there 
is a canonical isomorphism of $(\bar R\Cotimes\bar R)_{\sigma}$-modules
$$
\H^i_c(\Bunb_n^d,\; \Autb_{\bar E_1^*}^d
\Lotimes_{(\bar R\Cotimes\bar R)_{\sigma}}\Autb_{\bar E_2}^d)\,\iso\,
\left\{
\begin{array}{cc}
\bar R_{\sigma}, & \mbox{if}\; i=0\\
0, & \mbox{if} \; i\ne 0,
\end{array}
\right.
$$
where $\bar R_{\sigma}$ is viewed as a 
$(\bar R\Cotimes\bar R)_{\sigma}$-module
via the localized diagonal map 
$(\bar R\Cotimes\bar R)_{\sigma}\to\bar R_{\sigma}$.} 

\medskip\noindent
ii) The stack $\Bunb_n^d$ is not of finite type (for $n>1$). However, 
the cuspidality condition of Conjecture~\ref{Con_1} implies that
$\Autb_{\bar E}^d$ is the extension by zero from a substack of finite type
of $\Bunb_n^d$. So, in fact, we calculate the cohomology of a stack of finite 
type.

\medskip\noindent
iii) The definition of G. Laumon and L. Moret-Bailly (\cite{LMb}, 18.8) of the
cohomology with compact support of a stack is applicable here, because $\Bunb_n$
is a Bernstein-Lunts stack (cf. Remark~\ref{Rem_PGL}). 
\end{Rems}

\medskip\smallskip
\noindent
4.2 \ Essentially, the idea is to derive Main Global Theorem from Main 
Local Theorem (\cite{Ly2}). Actually, instead of using Main Local Theorem,
we will replace it by the following statement, 
which is easier to prove as soon as Conjecture~\ref{Con_1} is assumed true. 

 Let $\eta: \Bun^d_n\to\Pic^d X$ be the map that sends $L$ to 
$(\det L)\otimes\Omega^{(1-n)+(2-n)+\ldots+(n-n)}$. Let us write 
$\gamma: X^{(d)}\to\Pic^d X$ for the map that sends $D\in X^{(d)}$ to $\cO_X(D)$. 

\begin{Pp}
\label{Pp_summer}
 Assume that Conjecture~\ref{Con_1} is true. Then for each $d\ge 0$
there is a canonical isomorphism in $\D^b_c(\Pic^d X, \,
\bar R\Cotimes\bar R)_{\sigma}$
$$
\eta_!(\Aut^d_{\bar E_1^*}\otimes \varrho_! ({_n\cK^d_{\bar E_2}}))\,
\iso\, \gamma_!\HOM(\bar E_1,\bar E_2)^{(d)})[d+n^2(g-1)](\frac{d+n^2(g-1)}{2})
$$
\end{Pp}

 The complex $_n\cK^0_{\bar E}$ does not depend on $\bar E$ and will be 
denoted $_n\cK^0$. To prove the above proposition, we will only use 
the particular case $d=0$ of Main Local Theorem under the following form.

\begin{Lm} 
\label{Lm_MLTh_0}
There is a canonical isomorphism 
$ \eta_!\varrho_!({_n\cK^0}\otimes{_n\cK^0})\iso\gamma_!\bar R$ in 
$\D^b_c(\Pic^0 X, \bar R)$.
\end{Lm}
\begin{Prf} We have $\RG_c(\A^1,\cL_{\psi})=0$ not only in $\D^b_c(\Spec k,\Qlb)$
but also in $\D^b_c(\Spec k, \bar R)$.
This can be seen, for example, from (\cite{Ngo},
Lemma~3.3). Therefore, the proof given in (\cite{Ly2}, Lemma 6) holds 
for $\Qlb$-sheaves replaced by $\bar R$-sheaves.
\end{Prf}
 
\medskip 
  
\begin{Prf}\select{of Proposition~\ref{Pp_summer}}.
By Proposition~\ref{Pp_majus'}, 
$\DD(\Aut_{\bar E})\iso\Aut_{\bar E^*}$. Therefore,
$$
\eta_!(\Aut^d_{\bar E_1^*}\otimes \varrho_! ({_n\cK^d_{\bar E_2}}))\,\iso\,
\eta_!\DD\R\HOM(\varrho_!(_n\cK^d_{\bar E_2}), \Aut^d_{\bar E_1})
$$
Recall the map $\supp: {_n\Mod_d}\to X^{(d)}$  (cf. Sect.~2.2). 
Let $\pi:{_n\Mod_d}\to \Sh^d_0$ be the map that sends $(L\subset L')$ to $L'/L$. 
Consider the diagram 
$$
\begin{array}{ccccc}
&& _n\Mod_d \\
& \swarrow\lefteqn{\scriptstyle h^{\gets}} && 
\searrow\lefteqn{\scriptstyle h^{\to}}\\
\Bun_n &&&& \Bun_n,
\end{array}
$$ 
where $h^{\gets}$ (resp., $h^{\to}$) sends $(L\subset L')\in{_n\Mod_d}$ to $L$
(resp., to $L'$). Following \cite{FGV2}, consider the averaging functor
$
\H^d_{n,\bar E}: \D(\Bun_n)\to\D(\Bun_n)
$
given by
$$
\H^d_{n,\bar E}(K)=h^{\to}_!(h^{\gets *}K\otimes \pi^*\cL^d_{\bar
E})[dn](\frac{dn}{2})
$$
We have $\H^d_{n,\bar E}(\varrho_!({_n\cK^0}))\iso \varrho_!(_n\cK^d_{\bar E})$.
Define the functor 
$\tilde{\H}^{-d}_{n,\bar E}: \D(\Bun_n)\to
\D(X^{(d)}\times\Bun_n)$ by
$$
\tilde{\H}^{-d}_{n,\bar E}(K')=(\supp\times h^{\gets})_!(\pi^*\cL^d_{\bar E}
\otimes h^{\to *}K')[dn](\frac{dn}{2})
$$

 Let $\epsilon: X^{(d)}\times \Pic^0 X\to X^{(d)}\times \Pic^d X$ be the
isomorphism that sends $(D\in X^{(d)},\cA\in\Pic^0 X)$ to $(D, \cA(D))$. 
Denote by $p: X^{(d)}\times\Bun_n^0\to \Bun_n^0$ and 
$p': X^{(d)}\times\Pic^d X\to \Pic^d X$ the projections.
The next result is a straightforward application of the
formalism of six functors. 

\begin{Lm} For $K\in\D(\Bun^0_n)$ and $K'\in\D(\Bun^d_n)$ we have
$$
\eta_*\R\HOM(\H^d_{n,\bar E}(K), K')\,\iso\, p'_*\epsilon_*(\id\times\eta)_*
\R\HOM(p^*K, \tilde{\H}^{-d}_{n,\bar E^*}(K'))\eqno{\square}
$$ 
\end{Lm}

 Applying this lemma we get
$$
\eta_*\R\HOM(\varrho_!(_n\cK^d_{\bar E_2}), \Aut^d_{\bar E_1})\,\iso\,
p'_*\epsilon_*(\id\times\eta)_*\R\HOM(p^*\varrho_!(_n\cK^0),\,
\tilde{\H}^{-d}_{n,\bar E^*_2}(\Aut^d_{\bar E_1}))
$$ 
By (\cite{FGV2}, 9.5), we have
$$
\tilde{\H}^{-d}_{n,\bar E_2^*}(\Aut^d_{\bar E_1})\,\iso\,
(\bar E_1\otimes\bar
E_2^*)^{(d)}
\boxtimes \Aut^0_{\bar E_1} [d](\frac{d}{2})
$$
This yields an isomorphism
$$
\R\HOM(p^*\varrho_!(_n\cK^0),\,
\tilde{\H}^{-d}_{n,\bar E^*_2}(\Aut^d_{\bar E_1}))\,\iso\,
(\bar E_1\otimes\bar
E_2^*)^{(d)}\boxtimes \R\HOM(\varrho_!(_n\cK^0), \Aut^0_{\bar E_1}) 
[d](\frac{d}{2})
$$
Applying Lemma~\ref{Lm_MLTh_0}, we get
\begin{multline*}
\eta_*\R\HOM(\varrho_!(_n\cK^0), \Aut^0_{\bar E_1})\,\iso\;
\eta_*\DD(\varrho_!(_n\cK^0)\otimes \Aut^0_{\bar E_1^*})\,\iso\;
\DD\eta_!(\varrho_!(_n\cK^0)\otimes \Aut^0_{\bar E_1^*})\,\iso\\
(\DD\eta_!\varrho_!({_n\cK^0}\otimes{_n\cK^0}))[n^2(1-g)](\frac{n^2(1-g)}{2})\,
\iso\;\gamma_*\bar R[n^2(1-g)](\frac{n^2(1-g)}{2})
\end{multline*}
Since the diagram commutes
$$
\begin{array}{ccc}
X^{(d)}\times X^{(0)}&\iso & X^{(d)}\\
\downarrow\lefteqn{\scriptstyle \id\times\gamma} && 
\downarrow\lefteqn{\scriptstyle \gamma}\\
X^{(d)}\times\Pic^0 X & \toup{\epsilon} X^{(d)}\times\Pic^d X\toup{p'} & \Pic^d X,
\end{array}
$$
we get 
\begin{multline*}
\eta_*\R\HOM(\varrho_!(_n\cK^d_{\bar E_2}), \Aut^d_{\bar E_1})\,\iso\,
p'_*\epsilon_*((\bar E_1\otimes\bar E_2^*)^{(d)}\boxtimes\gamma_*\bar R)
[d+n^2(1-g)](\frac{d+n^2(1-g)}{2})\,\iso\\
\gamma_*(\bar E_1\otimes\bar E_2^*)^{(d)}[d+n^2(1-g)](\frac{d+n^2(1-g)}{2})
\end{multline*}
Applying the Verdier duality functor $\DD$, we get the desired assertion. \\
\end{Prf}(Proposition~\ref{Pp_summer})

\medskip\smallskip
\noindent
4.3 \ Pick a closed point $x\in X$. Recall that it defines the map $\alpha_x:\uPic^d
X\to \Pic^d X$ (cf. Sect.~3.4). Define the stacks 
$^x\Bun_n^d$ and $^x_n\cM_d$ from the cartesian squares
$$
\begin{array}{ccc}
^x_n\cM_d & \to & {_n\cM_d}\\
\downarrow\lefteqn{\scriptstyle \varrho_x} && 
\downarrow\lefteqn{\scriptstyle \varrho}\\
^x\Bun_n^d &\toup{\beta_x} & \Bun_n^d\\
\downarrow && \downarrow\lefteqn{\scriptstyle \eta}\\
\uPic^d X & \toup{\alpha_x} & \Pic^d X
\end{array}
$$
Note that the composition 
$^x\Bun_n^d\toup{\beta_x} \Bun_n^d\toup{\gr} \Bunb_n^d$ is a  
$\mu_n$-gerbe, so that $(\gr\comp\beta_x)_!\bar R\iso\bar R$. Thus, in Main Global
Theorem we may and will replace the calculation of cohomologies of $\Bunb_n^d$
by that of $^x\Bun^d_n$.

\begin{Def}
 Let $^x\Aut^d_{\bar E}$ be the pull-back of  
$\Autb^d_{\bar E}$ under $^x\Bun^d_n\to\Bunb_n^d$. Let also
$^x_n\cK^d_{\bar E}$ be the pull-back of $_n\cK^d_{\bar E}[1](\frac{1}{2})$
under $^x_n\cM_d \to  {_n\cM_d}$.  
\end{Def} 

 Proposition~\ref{Pp_summer} admits 
the following immediate corolary.

\begin{Cor}
\label{Cor_Pp_summer}
Assume that Conjecture~\ref{Con_1} is true. Then for any $d\ge 0$ there is a
canonical isomorphism in $\D^b_c(\Spec k, \,
\bar R\Cotimes\bar R)_{\sigma}$
\begin{multline*}
\;\;\;\RG_c({^x\Bun_n^d},\, {^x\Aut^d_{\bar E_1^*}}\otimes \varrho_{x!}(^x_n\cK^d_{\bar
E_2}))\iso\\
\!\!\RG_c({^xX^{(d)}},\, {^x\HOM(\bar E_1,\bar E_2)^{(d)}})[d+n^2(g-1)+2]
(\frac{d+n^2(g-1)+2}{2})\;\;\qquad\square
\end{multline*}
\end{Cor}

 For $d\ge 0$ define the complex $\cN^d$ on $^x\Bun_n^d$ as $\varrho_{x!}(\bar
R\Cotimes\bar R)[2d-2n^2(g-1)](d-n^2(g-1))$. Combining Corolary~\ref{Cor_Pp_summer} 
with Proposition~\ref{Pp_coh_of_xHOM}, we get the following result. 

\begin{Cor}
\label{Cor_closer}
Assume that Conjecture~\ref{Con_1} is true. Then for any $d>0$ there is a
canonical isomorphism in $\D^b_c(\Spec k, \,
\bar R\Cotimes\bar R)_{\sigma}$
$$
\tau_{\ge 1-d}\,\RG_c({^x\Bun_n^d},\, {^x\Aut^d_{\bar E_1^*}}
\otimes{^x\Aut^d_{\bar E_2}}\otimes \cN^d)\,\iso\, \bar R,
$$
where the $\bar
R\Cotimes\bar R$-module structure on
$\bar R$ is given via the diagonal map $\bar
R\Cotimes\bar R\to\bar R$. \QED
\end{Cor}

\begin{Prf}\select{of Main Global Theorem}\\
Recall that in \cite{FGV2} a vector bundle $M$ is called \select{very unstable}
if it can be represented as a direct sum of two vector bundles $M\iso M_1\oplus M_2$
with $\Ext^1(M_1,M_2)=0$. Let $\Bun_n^{vuns}\subset\Bun_n$ denote the substack of very
unstable vector bundles. By the cuspidality property, $*$-restriction of $\Aut_{\bar E}$ to
$\Bun_n^{vuns}$ vanishes. 

 Recall that we have fixed a line bundle $\cL^{est}$ on $X$ (cf. Sect.~2.4). There is
a constant $c_{g,n}$ such that for $d\ge c_{g,n}$ and $M\in \Bun_n^d$ the condition
$\Hom(M,\cL^{est})\ne 0$ implies that $M$ is very unstable. 

 Pick $c'\in\ZZ$ such that for each $d\in\ZZ$ all the nontrivial cohomology sheaves
of $\Autb^d_{\bar E}$ with respect to the usual $t$-structure are places in degrees $\le
c'$, the existence of such constant follows from formulae (\ref{iso_multb}). 

 Given a pair of integers $i$ and $d$, calculate $\H^i_c(\Bunb_n^d,\, 
\Autb^d_{\bar E_1^*}\otimes\Autb^d_{\bar E_2})$ as follows. Pick $k\in\ZZ$ large
enough, so that if we put $d'=d+kn$ then the following conditions are satisfied:
\begin{itemize}
\item[A)] $\;d'>0,\;\; d'\ge c_{g,n},\;\; d'>n^2(g-1)$
\item[B)] $\;i\ge 1-d',\;\; i>-2d'+2n^2(g-1)+2c'+2\dim({^x\Bun^{d'}_n})$
\end{itemize}
We have used the fact that the dimension of $^x\Bun^{d'}_n$ does not depend on $d'$. 

 Consider the map $\multb_{kx}:\Bunb_n^d\to\Bunb_n^{d'}$ that sends $L$ to $L(kx)$. By
Proposition~\ref{Pp_majus'}, we have 
$$
\multb^*_{kx}\Autb^{d'}_{\bar E}\iso \Autb^d_{\bar E}\boxtimes (\wedge^n \bar
E)^{\otimes k}_x
$$
Therefore, 
\begin{multline*}
\H^i_c(\Bunb_n^d,\, 
\Autb^d_{\bar E_1^*}\otimes\Autb^d_{\bar E_2})\otimes (\wedge^n \bar E_1^*
\otimes \wedge^n \bar E_2)^{\otimes k}_x\;\iso\;
\H^i_c(\Bunb_n^{d'},\, 
\Autb^{d'}_{\bar E_1^*}\otimes\Autb^{d'}_{\bar E_2})\;\iso\\
\H^i_c({^x\Bun_n^{d'}},\, 
{^x\Aut^{d'}_{\bar E_1^*}}\otimes{^x\Aut^{d'}_{\bar E_2}})
\end{multline*}
Recall the stack $\Bun_n^{est}$ defined in the proof of Proposition~\ref{Pp_majus'}. 
The condition A) implies that $\Aut^{d'}_{\bar E}$ is the extension by zero from 
$\Bun_n^{est}\cap \Bun_n^{d'}$. Let $U$ denote the preimage of 
$\Bun_n^{est}\cap \Bun_n^{d'}$ under $\beta_x: {^x\Bun^{d'}_n}\to \Bun_n^{d'}$. 
Then over $U$ the map $\varrho_x: {_n^x\cM_{d'}}\to{^x\Bun_n^{d'}}$ is a vector bundle
of rank $d'-n^2(g-1)$ with removed zero section. So, over $U$ we have
$$
\tau_{\ge -2d'+2n^2(g-1)+2} \;\cN^{d'}\,\iso\,  \bar R\Cotimes\bar R
$$
Now, using condition B), from Corolary~\ref{Cor_closer} we conclude that 
$$
\H^i_c({^x\Bun_n^{d'}}, {^x\Aut^{d'}_{\bar E_1^*}}
\otimes{^x\Aut^{d'}_{\bar E_2}})\;\iso\; \left\{
\begin{array}{cc}
\bar R_{\sigma}, & \mbox{if}\; i=0\\
0, & \mbox{if} \; i\ne 0,
\end{array}
\right.
$$
as $(\bar R\Cotimes\bar R)_{\sigma}$-modules. Since $(\wedge^n \bar E_1^*
\otimes \wedge^n \bar E_2)^{\otimes k}_x$ is a free $\bar R\Cotimes\bar R$-module
of rank 1, and 
$$
(\wedge^n \bar E_1^*
\otimes \wedge^n \bar E_2)^{\otimes k}_x\otimes_{\bar R\Cotimes\bar R} \bar R\;\iso\;
\bar R
$$
canonically, we are done. 
\end{Prf}

\appendix
\section{Some linear algebra}
\label{Some_linear_algebra}

In this appendix we prove Lemma~\ref{Lm_A.1} (cf. Sect. 3.4.1).
Let $A$ be a commutative ring of characteristic 0. Let $M$ be a
bounded complex of $A$-modules. Then $S_d$ acts on the complex
$\otimesl_{i=1}^d M$, and we put $\Sym^d(M)=(\otimesl_{i=1}^d M)^{S_d}$.
Clearly, $\Sym^{k+l}(M)$ is a direct summand of
$\Sym^k(M)\otimes_A \Sym^l(M)$, so that we have both natural maps
$$
\Sym^k(M)\otimes_A \Sym^l(M) \to \Sym^{k+l}(M)\;\;\mbox{ and }\;\;
\Sym^{k+l}(M) \to \Sym^k(M)\otimes_A \Sym^l(M)
$$
We will be interested in complexes $M$ of the form
$(\dots\to M^{-2}\to M^{-1}\to A\to 0)$, i.e., we suppose that $M^0=A$
and $M^i=0$ for $i>0$. Denote by $\sigma_{\le k},\sigma_{\ge k}$ the
'foolish' troncation functors. Put $V=\sigma_{\le -1} M$, so the sequence
of complexes 
$$
0\to A\to M\to V\to 0
$$
is exact. Define a morphism
$f_d:\Sym^d(M)\to \Sym^{d+1}(M)$ as the composition $A\otimes\Sym^d(M)\to
M\otimes\Sym^d(M)=\Sym^1(M)\otimes\Sym^d(M)\to \Sym^{d+1}(M)$.
We get an inductive system of complexes $(\Sym^d(M), f_d)_{d\in\NN}$.
Put $\Sym^{\infty}(M)=\varinjlim \Sym^d(M)$.
 
\begin{lm}
\label{B_Lm}
For any $d\ge 0$ the sequence of complexes
$$
0\to \Sym^d(M)\toup{f_d}\Sym^{d+1}(M)\to \Sym^{d+1}(V)\to 0
$$
is exact, where the second arrow is defined by functoriality
from the natural morphism $M\to V$.
\end{lm}
\begin{Prf}
If $k\le 0$ then $(\otimesl_{i=1}^{d+1} M)^k$ is the direct sum
$$
[\otimesl_{
\begin{array}{c}
\scriptstyle{i_1+\dots +i_{d+1}=k} \\ 
\scriptstyle{i_j=0\; \mbox{\scriptsize for some}\; j}
\end{array}
}
M^{i_1}\otimes\dots\otimes M^{i_{d+1}}] 
\oplus [(\otimesl_{i=1}^{d+1} V)^k]
$$
The group $S_{d+1}$ acts on every summand in square brackets.
It is easy to understand that the invariants
$$
[\otimesl_{
\begin{array}{c}
\scriptstyle{i_1+\dots +i_{d+1}=k} \\ 
\scriptstyle{i_j=0\; \mbox{\scriptsize for some}\; j}
\end{array}
}
M^{i_1}\otimes\dots\otimes M^{i_{d+1}}]^{S_{d+1}}
$$
are identified with $(\Sym^d(M))^k$.
\end{Prf}

\medskip

From the above lemma it follows that the natural map
$\Sym^d(M)\to \Sym^{\infty}(M)$ is injective and $\sigma_{\ge -d}
\Sym^d(M) \to \sigma_{\ge -d}\Sym^{\infty}(M)$ is an isomorphism for
$d\ge 0$. So, $\Sym^{\infty}(M)$ is a filtered complex with the
filtration $(\Sym^d(M))_{d\in \NN}$. Since the morphisms 
$\Sym^k(M)\otimes\Sym^l(M)\to \Sym^{k+l}(M)$ are compatible with
$f_d$, by passing to the limit we get the morphism of multiplication
$\Sym^{\infty}(M)\otimes \Sym^{\infty}(M)\to \Sym^{\infty}(M)$
(compatible with the above filtration).

\begin{lm}
Suppose that $M=(M^{-1}\toup{\eth} A)$, i.e., $M^i=0$ for $i<-1$.
Then $\Sym^{\infty}(M)$ is the Koszul
complex for $M$. In addition, 
$$
\Sym^d(M)\to
\sigma_{\ge -d}\Sym^{\infty}(M)
$$
is an isomorphism for any $d\ge 0$.  
\QED
\end{lm}

Now we impose on $M$ the additional condition:  $M^i=0$ for $i<-2$
and $M^{-2}=A$, so $M=(A\to M^{-1}\to A)$. Then we have
\begin{lm} Let $d\ge 0$.\\
1) For any $k$ we have $(\Sym^d(M))^k=(\Sym^d(M))^{-2d-k}$ canonically.\\
2) If $0\le k\le d$ then $(\Sym^d(M))^{-k}=\wedge^k(M^{-1})\oplus\wedge^{k-2}
(M^{-1})\oplus\wedge^{k-4}(M^{-1})\oplus\dots$ 
\QED
\end{lm}

Set $W=\sigma_{\ge -1} M$. Pick $c_1\in A^*$ and consider the exact
sequence $0\to W\to M\to A[2]\to 0$, where $M\to A[2]$ is given by $M^{-2}
\toup{c_1} A$. 
Let us define a morphism $g_d: \Sym^{d+1}(M)\to
\Sym^d(M)[2]$ as the composition $\Sym^{d+1}(M)\to
\Sym^1(M)\otimes\Sym^d(M)=M\otimes\Sym^d(M)\to (A[2])\otimes\Sym^d(M)$.
The proof of the next result is analogous to that of Lemma \ref{B_Lm}.
\begin{lm}
\label{Lm_B_sequence2}
For any $d\ge 0$ the sequence
$$
0\to \Sym^{d+1}(W)\to \Sym^{d+1}(M)\toup{g_d} \Sym^d(M)[2]\to 0
$$
is exact, where the first arrow is defined by functoriality from
the natural map $W\to M$. 
\QED
\end{lm}
\begin{Rem}
i) The map $(\Sym^{d+1}(M))^k\toup{g_d}(\Sym^d(M))^{k+2}$ is
described as follows. 
For $-d-1\le k\le 0$ this is the morphism
$$
\wedge^{-k}(M^{-1})\oplus\wedge^{-k-2}(M^{-1})\oplus\ldots\to
\wedge^{-k-2}(M^{-1})\oplus\ldots
$$ 
that sends $\wedge^{-k}(M^{-1})$ to zero
and induces isomorphisms on the others direct summands. 
For $k<-d-1$ this is an isomorphism preserving the direct
sum decomoposition.\\ 
ii) Passing to the limit we get an exact sequence
$$
0\to \Sym^{\infty}(W)\to \Sym^{\infty}(M) \toup{g} \Sym^{\infty}(M)[2]
\to 0
$$
\end{Rem}

Denote by $\xi_d$ the composition $\Sym^d(M)\hook{f_d}\Sym^{d+1}(M)
\toup{g_d} \Sym^d(M)[2]$. In particular, $\xi_1$ is given by the diagram
$$
\begin{array}{ccc}
& A &\!\! \to M^{-1}\to A\\
& \downarrow\lefteqn{\scriptstyle{c_1}}\\
A\to M^{-1}\to\!\! & A
\end{array}
$$
Define the morphism $\xi: \otimesl_{i=1}^d M\to \otimesl_{i=1}^d M[2]$
as $\xi_1\otimes\id\otimes\dots\otimes\id+\ldots+
\id\otimes\dots\otimes\id\otimes\xi_1$.
Then the diagram commutes
$$
\begin{array}{ccc}
\otimesl_{i=1}^d M & \toup{\xi} & \otimesl_{i=1}^d M[2]\\
\cup && \cup\\
\Sym^d(M) & \toup{\xi_d} & \Sym^d(M)[2]
\end{array}
$$

\begin{Prf}\select{of Lemma~\ref{Lm_A.1}}\\
$\Sym^{\infty}(W)$ is the Koszul complex for $\eth(e_1),\ldots,\eth(e_r)\in A$, so
that the natural map $\Sym^{\infty}(W)\to A/I$ is a quasi-isomorphism. Our
assertion follows now from Lemmas~\ref{Lm_B_sequence2} and \ref{B_Lm}.
\end{Prf}

\section{Chern classes}
\label{Chern_classes}

 Let $S$ be a smooth separated scheme of finite type, 
$\cA$ an invertible $\cO_S$-module, $f:Y\to S$ the total space of $\cA$
with removed zero section. Fix $n>0$ invertible as a function on 
$S$ (we don't need here $S$ to be defined over $k$).   

\begin{Lm_section}
\label{prop_C1}
The complex $f_!\mu_n$
is included into a distinguished triangle 
$f_!\mu_n\to \ZZ/n\ZZ[-2]
\toup{c(\cA)}
\mu_n$
on $S$, where $c(\cA)\in H^2(S,\mu_n)$ is the Chern class of 
$\cA$.
\end{Lm_section}

To prove this we need two lemmas.
\begin{Lm_section}
\label{C_Lm_triangle}
Let $\cA$ be an abelian category, $\D(\cA)$ be its derived category. Let
$K'\toup{\alpha} K\to K''$ be a distinguished triangle in $\D(\cA)$. 
Suppose that the morphism $\H^i(K'')\to \H^{i+1}(K')$ is surjective.
Then there is a unique morphism $\tau_{\le i+1}K'\to \tau_{\le i}K$ such 
that the composition $\tau_{\le i+1}K'\to \tau_{\le i}K\to\tau_{\le i+1}K$
is obtained from $\alpha$ by applying the functor $\tau_{\le i+1}$, and
the triangle $\tau_{\le i+1}K'\to \tau_{\le i}K\to\tau_{\le i}K''$
is distinguished.
\QED
\end{Lm_section}

\begin{Lm_section} 
\label{C_Lm}
The complex $f_*\mu_n$ is included into a
distinguished triangle 
$
f_*\mu_n \to \ZZ/n\ZZ[-1] \toup{c(\cA)}
\mu_n[1]
$,
where $c(\cA)\in \H^2(S,\mu_n)$ is the Chern class of $\cA$.
\end{Lm_section}
\begin{Prf}
The Kummer exact sequence $1\to \mu_n\to \Gm\toup{x\mapsto x^n}
\Gm\to 1$ on $S$ defines a distinguished triangle $\Gm[1]\to
\Gm[1]\toup{\delta}
\mu_n[2]$, and the Chern class of $\cA\in\Pic S=\Hom_{\D}(\ZZ,\Gm[1])$
is the composition $\delta\comp\cA\in \Hom_{\D}(\ZZ,\mu_n[2])$ of
morphisms in the derived category.

 Consider now the Kummer exact sequence on $Y$. It provides a
distinguished triangle $f_*\mu_n\to f_*\Gm\to f_*\Gm$ on $S$. As is easily 
seen, the moprhism $(\R^0f_*)\Gm\to (\R^1f_*)\mu_n$ is surjective (the
question is local for the \'etale topology on $S$, and one can assume $f$
to be the projection $S\times\Gm\to S$).
Since $\tau_{\le 1}f_*\mu_n \iso f_*\mu_n$, by Lemma~\ref{C_Lm_triangle} we
get a distinguished triangle
$$
 f_*\mu_n\to (\R^0f_*)\Gm \toup{\varkappa} (\R^0f_*)\Gm
$$
 Let us now construct a morphism $(\R^0f_*)\Gm\to \ZZ$. 
If $U$ is a
smooth scheme then 
$$
\H^0(U\times \Gm,\Gm)\iso \H^0(U,\Gm)\times
\H^0(U,\ZZ)
$$
canonically. Let $S'\to S$ be an \'etale morphism. Put
$Y'=S'\times_S Y$ and denote by $g:\Gm\times Y'\to Y'$ the action
of $\Gm$ on $Y'$. Let $s\in \H^0(S', (\R^0f_*)\Gm)$. Since $Y'$ 
is smooth, to $s\comp g$ there corresponds an element of
$\H^0(Y',\ZZ)=\H^0(S',\ZZ)$. This provides a morphism
$(\R^0f_*)\Gm\to \ZZ$. It is included into an exact sequence
$$
0\to \Gm\to (\R^0f_*)\Gm\to \ZZ\to 0,
$$ 
where the first arrow comes from the natural morphism $\Gm\to
f_*f^*\Gm$. Using Chech coverings one proves that the corresponding
element of $\Ext^1_S(\ZZ,\Gm)=\Pic S$ is $\cA$. In other words,
we get a distinguished triangle $(\R^0f_*)\Gm\to \ZZ\toup{\cA}
\Gm[1]$ on $S$. 

 The morphism $\varkappa$ yeilds a morphism of exact sequences
$$
\begin{array}{ccccccccc}
0& \to & \Gm & \to & (\R^0f_*)\Gm & \to & \ZZ & \to & 0\\
 &     &\downarrow\lefteqn{x\mapsto x^n} &&
\downarrow\lefteqn{\varkappa} && \downarrow\lefteqn{n} &&\\
0& \to & \Gm & \to & (\R^0f_*)\Gm & \to & \ZZ & \to & 0
\end{array}
$$
The latter provides a commutative diagram, where the rows and 
columns are distinguished triangles
$$
\begin{array}{ccccc}
(\R^0f_*)\Gm & \to & \ZZ & \toup{\cA} & \Gm[1]\\
\downarrow\lefteqn{\varkappa} && \downarrow\lefteqn{n} &&\downarrow\\
(\R^0f_*)\Gm & \to & \ZZ & \toup{\cA} & \Gm[1]\\
\downarrow && \downarrow && \downarrow\lefteqn{\delta}\\
f_*\mu_n[1] & \to & \ZZ/n\ZZ & \to & \mu_n[2]
\end{array}
$$
So, the morphism $\ZZ/n\ZZ  \to  \mu_n[2]$ in the lowest row is
$c(\cA)$.
\end{Prf}

\smallskip

Lemma~\ref{prop_C1} follows from Lemma~\ref{C_Lm} by Verdier duality.

\smallskip

{\small \scshape Universit\'e Paris-Sud, b\^at. 425, Math\'ematiques, 91405
Orsay France}

e-mail: {\tt Sergey.Lysenko@math.u-psud.fr}

\end{document}